\documentclass[12pt]{article}

\usepackage{amsfonts}
\usepackage{amssymb}
\usepackage{latexsym}
\usepackage{amsmath}
\usepackage{bibmods}

\usepackage{graphicx}
\usepackage{color}

\usepackage{amssymb}
\usepackage{amscd} 
\usepackage{color}

\newtheorem{theorem}{Theorem}[section] 
\newtheorem{corollary}[theorem]{Corollary}
\newtheorem{proposition}[theorem]{Proposition}

 \newtheorem{definition}{Definition}
 \newtheorem{remark}{Remark}

\newcommand{\rbx}{\hfill{\rule{1ex}{1ex}}}

\newcommand{\co}{{}{\scriptstyle{\mathcal O}}}

\begin{document}

\vspace*{10mm}

\begin{center}
{\Large\textbf{On Mellin convolution operators in Bessel potential
spaces}\footnote{This work was carried out when the second author
visited the Universiti Brunei Darussalam (UBD). The support of UBD
provided via Grant UBD/GSR/S\&T/19 is highly appreciated. \\ \indent
The work is also supported by the Georgian National Science
Foundation, Contract No. 31/39).}}
\end{center}

\vspace{5mm}

\begin{center}

\textbf{V. D. Didenko, R. Duduchava}

\vspace{2mm}

Universiti Brunei Darussalam, Bandar Seri Begawan, BE1410  Brunei;
\texttt{diviol@gmail.com}

I. Javakhishvili Tbilisi State University, Andrea Razmadze
Mathematical Institute,  University str. 2, Tbilisi 0186, Georgia;
\texttt{dudu@rmi.ge}

 \end{center}

  \vspace{10mm}




\vspace{10mm}

\begin{abstract}
Mellin convolution equations acting in Bessel potential spaces are
considered. The study is based upon two results. The first one
concerns the interaction of Mellin convolutions and Bessel potential
operators (BPOs). In contrast to the Fourier convolutions, BPOs and
Mellin convolutions do not commute and we derive an explicit formula
for the corresponding commutator in the case of Mellin convolutions
with meromorphic symbols. These results are used in the lifting of
the Mellin convolution operators acting on Bessel potential spaces
up to operators on Lebesgue spaces. The operators arising belong to
an algebra generated by Mellin and Fourier convolutions acting on
$\mathbb{L}_p$-spaces. Fredholm conditions and index formulae for
such operators have been obtained earlier by R. Duduchava and are
employed here. Note that the results of the present work find
numerous applications in boundary value problems for partial
differential equations, in particular, for equations in domains with
angular points.
 \end{abstract}

\section*{Introduction}

Boundary value problems for elliptic equations in domains with
angular points play an important role in applications and have a
rich and exciting history. A prominent representative of this family
is the Helmholtz equation. In the classical $\mathbb{W}^1$-setting,
the existence and uniqueness of the solution of coercive systems
with various type of boundary conditions are easily obtainable by
using the Lax-Milgram Lemma (see, e.g., \cite{DTT14} where
Laplace-Beltrami equations are considered on smooth surface with
Lipschitz boundary). Similar problems arise in new applications in
physics, mechanics and engineering. Thus recent publications on
nano-photonics \cite{BCC12a,GB10} deal with physical and engineering
problems described by BVPs for the Helmholts equation in $2D$
domains with angular points.  They are investigated with the help of
a modified Lax-Milgram Lemma for so called $T$-coercive operators.
Similar problems occur for the Lam\'e systems in elasticity,
Cauchy-Riemann systems, Carleman-Vekua systems in generalized
analytic function theory etc.

Despite an impressive number of publications and ever growing
interest to such problems, the results available to date are not
complete. In particular, serious difficulties arise if an
information on the solvability in non-classical setting in the
Sobolev spaces $\mathbb{W}^1_p$, $1<p<\infty$ is required, and one
wants to study the solvability of equivalent boundary integral
equations in the trace spaces $\mathbb{W}^{1-1/p}_p$ on the
boundary. Integral equations arising in this case often have fixed
singularities in the kernel and are of Mellin convolution type. For
example, \cite{BDKT13} describes how model BVP's in corners emerge
from the localization of BVP for the Helmholtz equation in domains
with Lipschitz boundary. Consequently, an attempt to study the
corresponding Mellin convolution operators in Bessel potential
spaces has been undertaken in \cite{Du13}. However, not all of the
results presented there are correct and one of the aims of this work
is to provide correct formulations and proofs. We also hope that the
results of the present paper will be helpful in further studies of
boundary value problems for various elliptic equations in Lipschitz
domains.

One such a model problem has been studied in \cite{DT13}. More
precisely, consider the following BVP with mixed Dirichlet--Neumann
boundary conditions,
\begin{equation}\label{eqn1}
\left\{\begin{array}{ll}
\Delta u(x) + k^2 u(x)=0, \quad & x\in\Omega_\alpha,\\[2mm]
u^+(t)=g(t),  & t\in\mathbb{R}^+,\\[2mm]
(\partial_{\boldsymbol{\nu}} u)^+(t)=h(t),  & t\in\mathbb{R}_\alpha,
\end{array}\right.
\end{equation}
in a corner $\Omega_\alpha$ of magnitude $\alpha$,
\begin{equation*}
\begin{aligned}
&\partial\Omega_\alpha=\mathbb{R}^+\cup\mathbb{R}_\alpha,\;\; \mathbb{R}^+=(0,\infty),\\[2mm]
&\mathbb{R}_\alpha:=\{te^{i\alpha}=(t\cos\,\alpha,t\sin\,\alpha)\;:\;
t\in\mathbb{R}^+\}.
\end{aligned}
\end{equation*}
By \cite{DT13} the BVP \eqref{eqn1} is reduced to the following
equivalent system of boundary integral equations on $\mathbb{R}^+$,
\begin{align}\label{eqn2}
\left\{\begin{array}{ll}\varphi+\displaystyle\frac1{2\pi}\left[\boldsymbol{K}^1_{e^{i\alpha}}
     +\boldsymbol{K}^1_{e^{-i\alpha}}\right]\psi=G_1,\\[3mm]
\psi-\displaystyle\frac1{2\pi}\left[\boldsymbol{K}^1_{e^{i\alpha}}
     +\boldsymbol{K}^1_{e^{-i\alpha}}\right]\varphi=H_1\quad      & \end{array}\right.
\end{align}
where
\begin{equation}\label{eqn3}
\boldsymbol{K}^1_{e^{\pm
i\alpha}}\psi(t):=\frac1\pi\int_0^\infty\frac{\psi(\tau)d\tau}{t
     -e^{\pm i\alpha}\tau},\qquad 0<|\alpha|<\pi,
 \end{equation}
are Mellin convolution operators with homogeneous kernels of order
$-1$ (see e.g. \cite{Du79,Du82} and Section \ref{s1} below), also
called integral equations with fixed singularities in the kernel.
Similar integral operators arise in the theory of singular integral
equations with the complex conjugation if the contour of integration
possess corner points. A complete theory of such equations is
presented in \cite{DL85,DLS95}, whereas various approximation
methods have been investigated in \cite{DiSi08,DiRS95,DiV96}. For a
more detailed survey of this theory, applications in elasticity, and
numerical methods for the corresponding equations we refer the
reader to
\cite{Du75a,Du77,Du79,Du82,Sc85,DH:2011a,DH:2013a,DS:2002}. Note
that a similar approach has been employed by M.~Costabel and
E.~Stephan \cite{Cos83,CS84} in order to study boundary integral
equations on curves with corner points.

Nevertheless, the results available are not sufficient in order to
solve the problems arising in the investigation of BVP \eqref{eqn1}.
First of all, we are looking for a solution to BVP \eqref{eqn3} in
the classical (finite energy) formulation
 \begin{equation}\label{eqn4}
\begin{gathered}
g\in\mathbb{H}^{1/2}(\mathbb{R}^+),\quad h\in\mathbb{H}^{-1/2}(\mathbb{R}_\alpha), \quad u\in\mathbb{H}^1(\Omega_\alpha)=\mathbb{W}^1(\Omega_\alpha),\\
u(x)=\co (1)\quad {\rm as}\quad|x|\to\infty \, ,
 \end{gathered}
 \end{equation}
or in the non-classical formulation
\begin{equation}\label{eqn5}
\begin{gathered}
g\in\mathbb{W}^{1-1/p}_p(\mathbb{R}^+),\quad h\in\mathbb{W}^{-1/p}_p(\mathbb{R}_\alpha), \quad u\in\mathbb{H}^1_p(\Omega_\alpha)=\mathbb{W}^1_p(\Omega_\alpha),\\
u(x)=\co(1)\quad {\rm as}\quad|x|\to\infty, \quad 1<p<\infty.
\end{gathered}
 \end{equation}
The non-classical formulation is very helpful to explore the maximal
smoothness of a solution to the BVP. This plays an important role in
approximation methods and other applications..

The corresponding equivalent system of boundary integral equation
\eqref{eqn2} must be considered in the Bessel potential space
$\widetilde{\mathbb{H}}^{-1/2}(\mathbb{R}^+)$ in the case of
classical setting \eqref{eqn4} or in the Besov (Sobolev-Slobodeckii)
space $\widetilde{\mathbb{W}}^{-1/p}(\mathbb{R}^+)$ in the case of
non-classical setting \eqref{eqn5}. While doing so one encounters
the three major tasks.
 \begin{itemize}
 \item
In general, Mellin convolution operators are not bounded in neither
Besov nor Bessel potential spaces. Therefore, in order to study
equations \eqref{eqn2} in the spaces of interest, one has to find a
subclass of multipliers with the boundedness property.
 \item
If boundedness criteria for the operators associated with equation
\eqref{eqn2} are available, one can lift this equation from the
Besov or the Bessel potential space to a Lebesgue space.
 \item
The  lifted equations should be studied in the Lebesgue space.
 \end{itemize}

A suitable class of Mellin convolution operators bounded in the
Bessel potential spaces was proposed in \cite{Du13}. These are
Mellin convolutions with {\em admissible meromorphic kernels} (see
\eqref{eqn20} below). Having proved the boundedness result, one can
study convolution equations in Bessel potential spaces. In
particular, by lifting an equation with Mellin convolution operator
$\mathfrak{M}_a^0$ with the help of Bessel potential operators
$\mathbf{\Lambda}^s_+$ and $\mathbf{\Lambda}^{s-r}_-$, one obtains
an equation in $\mathbb{L}_p$-space with the operator
$\mathbf{\Lambda}^{s-r}_-\mathfrak{M}_a^0 \mathbf{\Lambda}^{-s}_+$.
However, the resulting operator
$\mathbf{\Lambda}^{s-r}_-\mathfrak{M}_a^0 \mathbf{\Lambda}^{-s}_+$
is neither Mellin nor Fourier convolution and in order to describe
its properties, one first has to study the commutators of Bessel
potential operators and Fourier convolutions with discontinuous
symbols. As was already mentioned, this problem has been considered
in \cite{Du13}, but not all of the results of that work are correct.
Therefore, in Section~\ref{s1} the commutator problem is discussed
once again, and Theorem \ref{t3.1}, Corollary \ref{c3.2} below
provide correct formulae for the corresponding commutators.

The lifted operator
$\mathbf{\Lambda}^{s-r}_-\mathfrak{M}_a\mathbf{\Lambda}^{-s}_+$
belongs to the Banach algebra generated by Mellin and Fourier
convolution operators with discontinuous symbols. Such algebras have
been studied before in \cite{Du87} and the results obtained are
systematized and updated in the recent paper \cite{Du13}. In \S\,2,
these results are applied to the lifted equation, hereby
establishing properties of the initial Mellin convolution equation
in the Bessel potential space.

The results of the present paper are applied to BVPs for the
Helmholtz and Lam\'e equations in domains with corners and the
corresponding paper of R. Duduchava, M. Tsaava and T. Tsutsunava
will appear soon. These problems were investigated earlier only by
means of Lax-Milgram Lemma \cite{BCC12a}. In contrast to that, the
approach of the present work is more fruitful and provides better
tools to analyze the solvability of the equations involved and the
asymptotic behaviour of their solutions. Moreover, it can also be
used to study the Schr\"odinger operator on combinatorial and
quantum graphs. Such a problem has attracted a lot of attention
recently, since the operator mentioned has a wide range of
applications in nano-structures \cite{Ku04,Ku05} and possesses
interesting properties. Another area where the results of the
present paper can be useful, is the study of Mellin
pseudodifferential operators on graphs. This problem has been
considered in \cite{RR12} but in the periodic case only. Moreover,
some of the result obtained play an important role in the theory of
approximation methods for Mellin operators in Bessel potential
spaces.

The present paper is organized as follows. In the first two sections
we define Mellin convolution operators and recall some of their
properties. In the second section we also consider Fourier
convolution operators in the Bessel potential spaces and discuss the
lifting of these operators from the Bessel potential spaces to
Lebesgue spaces, mostly according the papers \cite{Du79,Es81}. For
Mellin convolutions such a lifting operation has not been studied
before, and in Section~$3$ the interaction between Bessel potential
operators and the Mellin convolution $\mathbf{K}_c^1$ with the
kernel $(t-c\tau)^{-1}$ is considered. In particular, we derive
formulae for commutators of Bessel potential operators and Mellin
convolutions, and these results are crucial for our further
considerations.

Section~$4$ recalls results from \cite{Du87,Du13} concerning the
Banach algebra generated by Fourier and Mellin convolution operators
in Lebesgue spaces with weight. These results, together with Theorem
\ref{t3.1} and Corollary \ref{c3.2}, are used in Section~$5$ in
order to describe the lifting of Mellin convolution operators from
the Bessel potential spaces up to operators in Lebesgue spaces. It
turns out that the objects arising belong to a Banach algebra
generated  by Mellin and Fourier convolutions in
$\mathbb{L}_p$-space on the semi-axis. The main result here is
represented by Theorem~\ref{t5.1} and Theorem~\ref{t5.2}, where the
interaction between Bessel potential operators and the Mellin
convolution resulting from the lifting of a model operator
$\mathbf{K}_c^1$ is described. Theorem \ref{t5.3} deals with the
lifting of the operator $\mathbf{K}_c^2$. In conclusion of
Section~$5$, we present explicit formulae for the symbols of Mellin
convolution operators with meromorphic kernels, which allow us to
find Fredholm criteria and an index formula for the operators under
consideration (see Theorem~\ref{t5.4} and Corollary~\ref{c5.5}).

 %
\section{Mellin convolution operators}
\label{s1}

Equations \eqref{eqn2} are a particular case of the Mellin
convolution equation
 \begin{equation}\label{eqn6}
\mathfrak{M}_a^0\mathbf{\varphi}(t):=c_0\mathbf{\varphi}(t)+\frac{c_1}{\pi
i}\int_0^\infty
     \frac{\mathbf{\varphi}(\tau)\,dt}{\tau-t}+\int_0^\infty\mathcal{K}\left(\frac t\tau\right)\mathbf{\varphi}(\tau)\frac{d\tau}{\tau}=f(t)
 \end{equation}
where $c_0,c_1\in\mathbb{C}$. If the kernel $\mathcal{K}$ satisfies
the condition
 \begin{equation*}
    \int_0^\infty t^\beta|\mathcal{K}(t)|\frac{dt}t<\infty, \quad 0<\beta<1,
 \end{equation*}
then both equation \eqref{eqn6} and analogous equations on the unit
interval $I:=(0,1)$  considered, respectively, on Lebesgue spaces
$\mathbb{L}_p(\mathbb{R}^+)$ and $\mathbb{L}_p(I)$, are fully
studied in \cite{Du79}.

Let $a$ be an essentially bounded measurable $N\times N$ matrix
function $a\in\mathbb{L}_\infty(\mathbb{R})$, and let
$\mathcal{M}_\beta$ and $\mathcal{M}^{-1}_\beta$ denote,
respectively, the Mellin transform and its inverse, i.e.
 \begin{align*}
\mathcal{M}_\beta\psi(\xi) & :=\int\limits_0^\infty
t^{\beta-i\xi}\psi(t)\,
    \frac{dt}t,\;\;\xi\in\mathbb{R},\\
\mathcal{M}^{-1}_\beta\varphi(t) &
:=\frac1{2\pi}\int\limits_{-\infty}^{\infty}
    t^{i\xi-\beta}\varphi(\xi)\,d\xi,\;\;t\in\mathbb{R}^+.
\end{align*}
On the Schwartz space  $\mathbb{S}(\mathbb{R}^+)$ of the fast
decaying functions on $\mathbb{R}^+$, consider the following
equation
 \begin{equation}\label{eqn7}
 \mathfrak{M}_a^0\varphi(t)=f(t),
 \end{equation}
where $\mathfrak{M}_a^{0}$ is the Mellin convolution operator,
  \begin{equation}\label{eqn8}
  \begin{aligned}
\mathfrak{M}_a^0\varphi(t):&=\mathcal{M}^{-1}_\beta
a\mathcal{M}_\beta \varphi(t)\\
&= \frac1{2\pi}\int\limits_{-\infty}^{\infty}
\!a(\xi)\int\limits_0^\infty\!
   \Big(\frac t\tau\Big)^{i\xi-\beta}\varphi(\tau)\,\frac{d\tau}\tau\,d\xi,\quad
 \varphi\in\mathbb{S}(\mathbb{R}^+).
\end{aligned}
 \end{equation}
Note that equation \eqref{eqn6} has the form \eqref{eqn7} with the
function $a$ defined by
 \[
 a(\xi):=c_0+c_1\coth\,\pi(i\beta+\xi)+(\mathcal{M}_\beta\mathcal{K})(\xi).
 \]

Equations of the form \eqref{eqn6}, \eqref{eqn7} and similar
equations on finite intervals often arise in various areas of
mathematics and mechanics (see \cite{Du79,Ka73}).

The function $a(\xi)$ in \eqref{eqn8} is usually referred to as the
symbol of the Mellin operator $\mathfrak{M}_a^{0}$. Further, if the
corresponding Mellin convolution operator $\mathfrak{M}_a^0$ is
bounded on the weighted Lebesgue space
$\mathbb{L}_p(\mathbb{R}^+,t^\gamma)$ endowed with the norm
 \[
\big\|\varphi\mid\mathbb{L}_p(\mathbb{R}^+,t^\gamma)\big\|:=\bigg[\int\limits_0^\infty
t^\gamma|\varphi(t)|^p\,dt\bigg]^{1/p},
 \]
then the symbol $a(\xi)$ is called a Mellin
$\mathbb{L}_{p,\gamma}$--multiplier.

The two most important examples of Mellin convolution operators are
\begin{equation*}
S_{\mathbb{R}^+}\varphi(t):=\frac1{\pi i}
\int\limits_0^\infty\frac{\varphi(\tau)\,d\tau}{\tau-t}, \qquad
\boldsymbol{K}^m_c\varphi(t):=\frac1{\pi
i}\int\limits_0^\infty\frac{\tau^{m-1}\varphi(\tau)\,d\tau}{(t
     -c\,\tau)^m},
 \end{equation*}
where ${\rm Im}\,c\not=0$ and $m\in \mathbb{N}$ (see \eqref{eqn3},
\eqref{eqn6}). The operator $S_{\mathbb{R}^+}$ is the celebrated
Cauchy singular integral operator. The Mellin symbols of these
operators are (cf. \cite[\S\, 2]{Du13})
\begin{equation*}
\begin{array}{l}
 \hskip-7mm \sigma(S_{\mathbb{R}^+})(\xi):=-i\cot\pi(\beta-i\xi),\quad \xi\in\mathbb{R},\\[2ex]
\hskip-7mm\sigma(\boldsymbol{K}^m_c)(\xi):=\displaystyle\binom{\beta-i\xi-1}{m-1}\,\displaystyle
  \frac{e^{\mp\pi(\beta-i\xi)i}i}{\sin\pi(\beta-i\xi)}\,c^{\beta-i\xi-m}, \quad0<\pm\arg c<\pi,
\end{array}
 \end{equation*}
 where
  $$
 \displaystyle\binom{\theta-1}{m-1}:=\displaystyle\frac{(\theta-1)\cdots(\theta-m+1)}{(m-1)!},
     \quad\binom{\theta-1}{0}:=1.
  $$
In particular,
\begin{align}\label{eqn9}
\mathcal{M}_\beta\mathcal{K}^1_{-c}(\xi)&
     =\frac{c\,^{\beta-i\xi-1}}{\sin\pi(\beta-i\xi)},\qquad -\pi<\arg c<\pi, \\
\label{eqn10} \mathcal{M}_\beta\mathcal{K}^1_{-1}(\xi) &
=\frac1{\sin\pi (\beta-i\xi)}, \qquad\xi\in\mathbb{R}.
\end{align}

The study of the equation \eqref{eqn7} does not require much effort.
The Mellin transform $\mathcal{M}_\beta$ converts \eqref{eqn7} into
the equation
  \begin{equation}\label{eqn11}
  a(\xi)\mathcal{M}_\beta\varphi(\xi)=\mathcal{M}_\beta f(\xi).
 \end{equation}
If $\inf|\det\,a(\xi)|>0$ and the matrix-function $a^{-1}$ is a
Mellin $\mathbb{L}_{p,\gamma}$-multiplier, then equation
\eqref{eqn11} has the unique solution
$\varphi=\mathcal{M}^0_{a^{-1}}f$.

The solvability of analogues of equation \eqref{eqn8} on the unit
interval $I=(0,1)$ in a weighted Lebesgue space
$\mathbb{L}_p([0,1],t^\gamma)$ is also well understood. Thus if
 \begin{align}\label{eqn12}
1<p<\infty, \qquad -1<\gamma<p-1, \qquad
\beta:=\displaystyle\frac{1+\gamma}p, \qquad
    0<\beta<1,
 \end{align}
then one can use the isomorphisms
 \begin{equation}\label{eqn13}
\begin{aligned}
Z_\beta&:\mathbb{L}_p([0,1],t^\gamma)\rightarrow\mathbb{L}_p(\mathbb{R}^+),
&
     Z_\beta\varphi(\xi):=e^{-\beta\xi}\varphi(e^{-\xi}),&\quad \xi\in\mathbb{R}^+,\\[1ex]
Z^{-1}_\beta&:\mathbb{L}_p(\mathbb{R}^+)\rightarrow\mathbb{L}_p([0,1],t^\gamma),&
     Z^{-1}_\beta\psi(t):=t^{-\beta}\psi(-\ln\,t), &\quad t\in(0,1],
\end{aligned}
\end{equation}
and transform the corresponding equation on the unit interval $I$
into an equivalent Wiener-Hopf equation, i.e. into the equation
 \begin{equation}\label{eqn14}
W_{\mathcal{A}_\beta}\psi(x)=c_0\psi(x)+\int\limits_0^\infty
\mathcal{K}_1(x-y)\varphi(y)dy
     =f_0(t).
 \end{equation}
The Fourier transform of the kernel $\mathcal{K}_1$ is called the
symbol of the corresponding Fourier convolution operator and is used
to describe Fredholm properties, index and solvability of the
equation \eqref{eqn14}. In passing note that Fourier convolution
equations with discontinuous symbols are well studied
\cite{Du75a,Du77,Du78,Du79,Th85}.

 %
\section{Fourier convolution operators in the Bessel potential spaces:
 definition and lifting}
\label{s2}

Let $N$ be a positive integer and let $\mathfrak{A}$ be a Banach
algebra. If no confusion can arise, we write $\mathfrak{A}$ for both
scalar and matrix $N\times N$ algebras with entries from
$\mathfrak{A}$. Similarly, the same notation $\mathfrak{A}$ is used
for the set of $N$-dimensional vectors with entries from
$\mathfrak{A}$.  It will be usually clear from the context what kind
of space or algebra is considered.

Along with Mellin convolutions $\mathfrak{M}_a^0$, let us consider
the Fourier convolution operators
$$
W^0_a\varphi:=\mathcal{F}^{-1}a\mathcal{F}\varphi, \;\;
\varphi\in\mathbb{S}(\mathbb{R}),
$$
where $a\in\mathbb{L}_{\infty,loc}(\mathbb{R})$ is a locally bounded
$N\times N$ matrix function, called the symbol of $W_a^0$ and
$\mathcal{F}$ and $\mathcal{F}^{-1}$ are, respectively, the direct
and inverse Fourier transforms, i.e.
 \[
 \mathcal{F}\varphi(\xi):=\int\limits_{-\infty}^\infty e^{i\xi x} \varphi(x)\,dx,\quad
\mathcal{F}^{-1}\psi(x):=\frac1{2\pi}\int\limits_{-\infty}^\infty
e^{-i\xi x}\psi(\xi)\,d\xi,\;\;x\in\mathbb{R}.
 \]

Let $1<p<\infty$. An $N\times N$ matrix symbol $a(\xi)$,
$\xi\in\mathbb{R}$ is called $\mathbb{L}_p$-multiplier if the
corresponding convolution operator
$W^0_a\;:\;\mathbb{L}_p(\mathbb{R})
\longrightarrow\mathbb{L}_p(\mathbb{R})$ is bounded. The set of all
$\mathbb{L}_p$-multipliers is denoted by
$\mathfrak{M}_p(\mathbb{R})$. It is known  (see, e.g. \cite{Du79}),
that $\mathfrak{M}_p(\mathbb{R})$ is a Banach subalgebra of
$\mathbb{L}_\infty(\mathbb{R})$ which contains the algebra
$\mathbf{V}_1(\mathbb{R})$ of all functions with finite variation.
For $p=2$ we have the exact equality
$\mathfrak{M}_2(\mathbb{R})=\mathbb{L}_\infty(\mathbb{R})$.

The operator
 \begin{eqnarray*}
W_a:=r_{\mathbb{R}^+}W^0_a\;:\;\mathbb{L}_p(\mathbb{R}^+)\longrightarrow
     \mathbb{L}_p(\mathbb{R}^+),
 \end{eqnarray*}
where
$r_{\mathbb{R}^+}:\mathbb{L}_p(\mathbb{R})\to\mathbb{L}_p(\mathbb{R}^+)$
denotes the restriction operator, is called the con\-volution on the
semi-axis $\mathbb{R}^+$ or the Wiener-Hopf operator. It is worth
noting that unlike the operators $W^0_a$ and $\mathfrak{M}^0_a$,
which possess the property
 \begin{align}\label{eqn15}
W^0_aW^0_b=W^0_{ab},\quad
\mathfrak{M}^0_a\mathfrak{M}^0_b=\mathfrak{M}^0_{ab}\quad \text{for
all}\quad a,b\in\mathfrak{M}_p(\mathbb{R}),
 \end{align}
the product of Wiener--Hopf operators cannot be computed by the
simple rule \eqref{eqn15}. Thus for the operators  $W_a$ and $W_b$,
a similar relation
\begin{equation}\label{eqn16}
W_aW_b=W_{ab}
\end{equation}
is valid if and only if either $a(\xi)$ has an analytic extension
into the lower half plane or $b(\xi)$ has an analytic extension into
the upper half plane \cite{Du79}.

If conditions \eqref{eqn12} hold, the isometrical isomorphisms
\eqref{eqn13} are extended to the following isomorphisms of Lebesgue
spaces
 \begin{equation*}
\begin{aligned}
Z_\beta&:\mathbb{L}_p(\mathbb{R}^+,t^\gamma)\rightarrow\mathbb{L}_p(\mathbb{R}),
&
     Z_\beta\varphi(\xi):=e^{-\beta\xi}\varphi(e^{-\xi}),&\quad \xi\in\mathbb{R},\\[1ex]
Z^{-1}_\beta&:\mathbb{L}_p(\mathbb{R})\rightarrow\mathbb{L}_p(\mathbb{R}^+,t^\gamma),&
     Z^{-1}_\beta\psi(t):=t^{-\beta}\psi(-\ln\,t), &\quad t\in\mathbb{R}^+,
\end{aligned}
\end{equation*}
and provide the following connection between the Fourier and Mellin
transformations and the corresponding convolution operators--viz.,
\begin{eqnarray*}
\begin{array}{c}
\mathcal{M}_\beta=\mathcal{F} {\bf
Z}_\beta,\qquad\mathcal{M}_\beta^{-1}
    ={\bf Z}^{-1}_\beta\mathcal{F}^{-1},\\[3mm]
\mathfrak{M}_a^{0}=\mathcal{M}_\beta^{-1} a\mathcal{M}_\beta={\bf
Z}^{-1}_\beta
    \mathcal{F}^{-1} a\mathcal{F} {\bf Z}_\beta={\bf Z}^{-1}_\beta W^0_a {\bf Z}_\beta.
\end{array}
\end{eqnarray*}
These identities also justify the following assertion.
 %
\begin{proposition}[\cite{Du79}]\label{p1}
Let $1<p<\infty$ and  $-1<\gamma<p-1$. The class of Mellin
$\mathbb{L}_{p,\gamma}$-multipliers does not depend on the parameter
$\gamma$ and coincides with the Banach algebra
$\mathfrak{M}_p(\mathbb{R})$ of Fourier $\mathbb{L}_p$-multipliers.
\end{proposition}
 %
\begin{corollary}[\cite{Du79}]\label{c1}
A Mellin convolution operator $\mathfrak{M}_a^0:
\mathbb{L}_p(\mathbb{R}^+, t^\gamma) \to \mathbb{L}_p(\mathbb{R}^+,
t^\gamma)$ of \eqref{eqn8} is bounded if and only if
$a\in\mathfrak{M}_p(\mathbb{R})$.
\end{corollary}

For $s\in\mathbb{R}$ and $1<p<\infty$, the  Bessel potential space,
known also as a fractional Sobolev space, is a subspace of the
Schwartz space $\mathbb{S}'(\mathbb{R})$ of the distributions having
the finite norm
 \begin{eqnarray*}
\|\varphi|\mathbb{H}^s_p(\mathbb{R})\|:=\left[\int_{-\infty}^\infty
     \left|\mathcal{F}^{-1}(1+|\xi|^2)^{s/2}
     (\mathcal{F}\varphi)(t)\right|^pdt\right]^{1/p}<\infty.
 \end{eqnarray*}

For the integer parameters $s=m\in \mathbb{N}$, space
$\mathbb{H}^s_p(\mathbb{R})$ coincides with the Sobolev space
$\mathbb{W}^m_p(\mathbb{R})$ endowed with an equivalent norm
 \begin{align*}
\|\varphi|\mathbb{W}^m_p(\mathbb{R})\|:=\left[\sum_{k=0}^m\int_{-\infty}^\infty
     \left|\frac{d^k\varphi(t)}{dt^k}\right|^pdt\right]^{1/p}.
 \end{align*}

If $s<0$, one gets the space of distributions. Moreover,
$\mathbb{H}^{-s}_{p'}(\mathbb{R})$ is the dual  to the space
$\mathbb{H}^s_p(\mathbb{R}^+)$, provided that
$p':=\displaystyle\frac p{p-1}$, $1<p<\infty$. Note that
$\mathbb{H}^s_2(\mathbb{R})$ is a Hilbert space with the inner
product
\begin{equation*}
    \langle\varphi, \psi\rangle_s =\int_{\mathbb{R}}
    (\mathcal{F}\varphi)(\xi) \overline{(\mathcal{F}\psi)(\xi)}
    (1+\xi^2)^s d\xi\;, \quad \varphi, \psi \in \mathbb{H}^s(\mathbb{R}).
\end{equation*}

By $r_\Sigma$ we denote the operator restricting functions or
distributions defined  on $\mathbb{R}$ to the subset
$\Sigma\subset\mathbb{R}$. Thus
$\mathbb{H}^s_p(\mathbb{R}^+)=r_{\mathbb{R}^+}(\mathbb{H}^s_p(\mathbb{R}))$,
and the norm in $\mathbb{H}^s_p(\mathbb{R}^+)$ is defined by
 \[
\|f|\mathbb{H}^s_p(\mathbb{R}^+\|=\inf_\ell \|\ell
f|\mathbb{H}^s_p(\mathbb{R})\|,
 \]
where $\ell f$ stands for any extension of $f$ to a distribution in
$\mathbb{H}^s_p(\mathbb{R})$.

Further, we denote by $\widetilde{\mathbb{H}}^s_p(\mathbb{R}^+)$ the
(closed) subspace of $\mathbb{H}^s_p(\mathbb{R})$ which consists of
all distributions supported in the closure of $\mathbb{R}^+$.

Note that $\widetilde{\mathbb{H}}^s_p(\mathbb{R}^+)$ is always
continuously embedded in $\mathbb{H}^s_p(\mathbb{R}^+)$ and for
$s\in(1/p-1,1/p)$ these two spaces coincide. Moreover,
$\mathbb{H}^s_p(\mathbb{R}^+)$ may be viewed as the quotient-space
$\mathbb{H}^s_p(\mathbb{R}^+):=\mathbb{H}^s_p(\mathbb{R})/\widetilde{\mathbb{H}}^s_p(
\mathbb{R}^-)$, $\mathbb{R}^-:=(-\infty,0)$.

If the Fourier convolution operator (FCO) on the semi-axis
$\mathbb{R}^+$ with the symbol
$a\in\mathbb{L}_{\infty,loc}(\mathbb{R})$ is bounded in the space
setting
 \begin{equation*}
W_a:=r_{\mathbb{R}^+}W^0_a\;:\;\widetilde{\mathbb{H}}^s_p(\mathbb{R}^+)\longrightarrow
     \mathbb{H}^{s-r}_p(\mathbb{R}^+).
 \end{equation*}
we say that $W_a$ has order $r$ and $a$ is an $\mathbb{L}_p$
multiplier of order $r$. The set of all $\mathbb{L}_p$ multipliers
of order $r$ is denoted by $\mathfrak{M}^r_p(\mathbb{R})$. Let us
mention another description of the space
$\mathfrak{M}^r_p(\mathbb{R})$, viz.
$a\in\mathfrak{M}^r_p(\mathbb{R})$ if and only if
$\lambda^{-r}a\in\mathfrak{M}_p(\mathbb{R})=\mathfrak{M}^0_p(\mathbb{R})$,
where $\lambda^r(\xi):=(1+|\xi|^2)^{r/2}$.

Note, that FCOs are particular cases of pseudodifferential operators
($\Psi$DOs).
 %
 \begin{theorem}\label{t1.1}
Let  $1<p<\infty$. Then
 \begin{enumerate}
 \item
For any $r,s\in\mathbb{R}$ and for any $\gamma\in \mathbb{C}$, ${\rm
Im}\,\gamma>0$,  pseudodifferential operators
$\mathbf{\Lambda}_{\gamma}^r:=\mathbf{\Lambda}_{+\gamma}^r$ and
$\mathbf{\Lambda}_{-\gamma}^r$ defined by
 \begin{equation}\label{eqn17}
  \begin{array}{l}
\mathbf{\Lambda}_\gamma^r=W_{\lambda^r_\gamma}\;:\;\widetilde{\mathbb{H}}^s_p(\mathbb{R}^+)
     \rightarrow\widetilde{\mathbb{H}}^{s-r}_p(\mathbb{R}^+),\\[3mm]
\mathbf{\Lambda}_{-\gamma}^r=W_{\lambda^r_{-\gamma}}\;:\;\mathbb{H}^s_p(\mathbb{R}^+)
     \rightarrow\mathbb{H}^{s-r}_p(\mathbb{R}^+),
\end{array}
 \end{equation}
where $\lambda^r_{\pm\gamma}(\xi):=(\xi\pm\gamma)^r$,
$\xi\in\mathbb{R}^+$, are isomorphisms between the corresponding
spaces.
\item
For any  operator
$\mathbf{A}\;:\;\widetilde{\mathbb{H}}^s_p(\mathbb{R}^+) \to
\mathbb{H}^{s-r}_p( \mathbb{R}^+)$ of order $r$, the following
diagram is commutative
 \begin{equation}\label{eqn18}
\begin{array}{ccc}\widetilde{\mathbb{H}}^s_p(\mathbb{R}^+) & \stackrel{\mathbf{A}}{\longrightarrow}
    &\mathbb{H}^{s-r}_p(\mathbb{R}^+)\\
\uparrow \mathbf{\Lambda}^{-s}_\gamma&&\downarrow \mathbf{\Lambda}_{-\gamma}^{s-r}\\
     \mathbb{L}_p(\mathbb{R}^+)&\stackrel{\mathbf{\Lambda}_{-\gamma}^{s-r}\mathbf{A}
     \mathbf{\Lambda}^{-s}_\gamma}{\longrightarrow}& \mathbb{L}_p(\mathbb{R}^+).\end{array}
 \end{equation}
Thus the diagram \eqref{eqn18} provides an equivalent lifting of the
operator $\mathbf{A}$ of order $r$ up to the operator
$\mathbf{\Lambda}_{-\gamma}^{s-r}\mathbf{A}
\mathbf{\Lambda}^{-s}_\gamma\;:\;\mathbb{L}_p(\mathbb{R}^+)\longrightarrow
\mathbb{L}_p(\mathbb{R}^+)$ of order $0$.
 \item Let $\mu, \nu\in \mathbb{R}$. If $a$ is an $\mathbb{L}_p$-multiplier
of order $r$, then for any complex numbers $\gamma_1, \gamma_2$ such
that ${\rm Im}\,\gamma_j>0$, $j=1,2$, the operator
$\mathbf{\Lambda}_{-\gamma_1}^\mu
W_a\mathbf{\Lambda}_{\gamma_2}^\nu$ is a Fourier convolution
$W_{a_{\mu,\nu}}$ of order $r+\mu+\nu$,
\begin{equation}\label{eqn19}
W_{a_{\mu,\nu}}\;:\;\widetilde{\mathbb{H}}^{s+\nu}_p(\mathbb{R}^+)
 \longrightarrow\mathbb{H}^{s-r-\mu}_p(\mathbb{R}^+),
\end{equation}
with the symbol
 $$
a_{\mu,\nu}(\xi):=(\xi-\gamma_1)^\mu a(\xi)(\xi+\gamma_2)^\nu .
 $$
In particular,  the lifting of the operator  $W_a$ up to the
operator $\mathbf{\Lambda}_{-\gamma}^{s-r}W_a
\mathbf{\Lambda}^{-s}_\gamma$ acting in the space
$\mathbb{L}_p(\mathbb{R}^+)$ is FCO of order zero with the symbol
\begin{equation*}
a_{s-r,-s}(\xi)=\lambda^{s-r}_{-\gamma}(\xi)a(\xi)\lambda^{-s}_\gamma(\xi)
    =\Big(\frac{\xi-\gamma}{\xi+\gamma}\Big)^{s-r}\,\frac{a(\xi)}{(\xi+\gamma)^r}\,.
\end{equation*}
\item
The Hilbert transform
$\mathbf{K}^1_1=iS_\mathbb{\mathbb{R}^+}=W_{{\rm -i\,sign}}$ is a
Fourier convolution operator and
\begin{gather*}
\mathbf{\Lambda}^s_{-\gamma_1}\boldsymbol{K}^1_1\mathbf{\Lambda}^{-s}_{\gamma_2}
     =W_{i\,g^s_{-\gamma_1,\gamma_2}\,{\rm sign}},
\end{gather*}
where
\begin{align*}
g^s_{-\gamma_1,\gamma_2}(\xi):=\left(\frac{\xi-\gamma_1}{\xi+\gamma_2}\right)^s.
\end{align*}
 \end{enumerate}
 \end{theorem}
\textbf{Proof.}   For the proof of items $(i)-(iii)$ we refer the
reader to \cite[Lemma 5.1]{Du79} and \cite{DS93,Es81}. The item
$(iv)$ is a consequence of $(ii)-(iii)$ (see \cite{Du79,Du13}). \rbx

Note that the operator equality in \eqref{eqn19} is in fact a
consequence of the relation \eqref{eqn16}.

 %
\section{Mellin convolution operators in the Bessel potential spaces--lifting}
\label{s3}

In contrast to the Fourier convolution operators the lifted Mellin
convolution operator is not a Mellin convolution anymore. Moreover,
there are Mellin convolution operators $\mathfrak{M}^0_{a_\beta}$
with symbols $a_\beta\in\mathfrak{M}_p(\mathbb{R})$ which are
unbounded in the Bessel potential spaces. Thus in order to study the
Mellin convolutions in the space of Bessel potentials, one has to
address the boundedness problem first. To this end, a class of
integral operators with admissible kernels was introduced in
\cite{Du13}. For the sake of simplicity, here we consider a lighter
version of such kernels.
 %
\begin{definition}
The function $\mathcal{K}$ is called an {\em admissible meromorphic
kernel} if  it can be represented in the form
\begin{equation}\label{eqn20}
\mathcal{K}(t):=\sum_{j=0}^\ell \frac{d_j}{t-c_j}+\sum_{j=\ell+1}^N
    \frac{d_j}{(t-c_j)^{m_j}},
 \end{equation}
where $d_j,c_j\in\mathbb{C}$, $j=0,1,\dots,N$,
$m_{\ell+1},\ldots,m_N \in \{2,3,\ldots\}$, and $0<\alpha_k:=|\arg
c_k|\leqslant \pi$ for $k=\ell+1,\dots,N$.
\end{definition}
Note that the kernel $\mathcal{K}(t)$ has poles at the points
$c_0,c_1,\ldots, c_N\in\mathbb{C}$.

Recall that boundary integral operators for BVPs in planar domains
with corners have admissible kernels (see \eqref{eqn2} and
\cite{Du79,Du82,Du13,DT13}).
 %
\begin{theorem}[{\cite[Theorem 2.5 and  Corollary 2.6]{Du13}}]\label{t1.2}
Let $1<p<\infty$ and $s\in\mathbb{R}$. If $\mathcal{K}$ is an
admissible kernel, then the Mellin convolution operator
\begin{equation}\label{eqn21}
\mathfrak{M}^0_a:\widetilde{\mathbb{H}}^s_p(\mathbb{R}^+)\longrightarrow
\mathbb{H}^s_p(\mathbb{R}^+),
 \end{equation}
where $a_\beta=\mathcal{M}_\beta \mathcal{K}$, is bounded.
\end{theorem}

The next result is crucial to what follows. Note that a similar
assertion appears in \cite{Du13}, but the proof contains fatal
errors.
 %
\begin{theorem}\label{t3.1}
Let $s\in\mathbb{R}$, $c,\gamma \in \mathbb{C}$,
$-\pi<\arg\,c\leqslant\pi$, $\arg\,c\not=0$, $0<\arg \gamma<\pi$ and
$-\pi<\arg(c\,\gamma)<0$. Then
\begin{gather}\label{eqn22}
\mathbf{\Lambda}^s_{-\gamma}{\bf K}^1_c=c^{-s}{\bf
K}^1_c\mathbf{\Lambda}^s_{-c\,\gamma},
\end{gather}
where $c^{-s}=|c|^{-s}e^{-s\arg (c)\,i}$.
\end{theorem}

\noindent {\em Proof of Theorem \ref{t3.1}.}  Taking into account
the mapping properties  of Bessel potential operators \eqref{eqn17}
and the mapping properties of a Mellin convolution operator with an
admissible kernel \eqref{eqn21}, one observes that both operators
\begin{align}\label{eqn23}
\begin{array}{lcr}
\mathbf{\Lambda}^s_{-\gamma}{\bf
K}^1_c&:&\widetilde{\mathbb{H}}^r_p(\mathbb{R}^+)
     \longrightarrow\mathbb{H}^{r-s}_p(\mathbb{R}^+),\\[1mm]
{\bf
K}^1_c\mathbf{\Lambda}^s_{-c\,\gamma}&:&\widetilde{\mathbb{H}}^r_p(\mathbb{R}^+)
     \longrightarrow\mathbb{H}^{r-s}_p(\mathbb{R}^+)
\end{array}
 \end{align}
are correctly defined and bounded for all $s\in\mathbb{R}$,
$1<p<\infty$, since  $0<\arg \gamma<\pi$ and
$0<-\arg(c\,\gamma)<\pi$.

On the other hand, let us note that the inverse superpositions ${\bf
K}^1_c\mathbf{\Lambda}^s_\gamma$ and
$\mathbf{\Lambda}^s_{c\,\gamma}{\bf K}^1_c$ are correctly defined
only for $1/p-1<s<1/p$ and $s=1,2,\ldots$.

For a smooth function with compact support  $\varphi\in
C^\infty_0(\mathbb{R}^+)$ and for $k=1,2,\ldots$ we can use
integration by parts and obtain
\begin{align}\label{eqn24}
    \frac{d^k}{dt^k}\,\boldsymbol{K}^1_{c}\varphi(t) & =\frac1\pi\int\limits_0^\infty \frac{d^k}{dt^k}\,\frac1{t-c\,\tau}\varphi(\tau)\,d\tau=\frac{(-c)^{-k}}\pi
    \int\limits_0^\infty \frac{d^k}{d\tau^k}\,\frac1{t-c\,\tau}\varphi(\tau)\,d\tau= \nonumber \\
&
=\frac{(-c)^{-k}}\pi\int\limits_0^\infty\frac1{t-c\,\tau}\,\frac{d^k\varphi(\tau)}{
    d\tau^k}\,d\tau=(-c)^{-k}\Big(\boldsymbol{K}^1_c\,\frac{d^k}{dt^k}\varphi\Big)(t).
\end{align}

Let us consider the case where $s$ is a positive integer, i.e.
$s=m=1,2,\ldots$. The Bessel potentials
$\mathbf{\Lambda}_\pm^m=W_{\lambda^m_{\pm\gamma}}$ are the Fourier
convolutions of order $m$ and they represent ordinary differential
operators of the order $m$, namely,
\begin{equation}\label{eqn25}
\mathbf{\Lambda}_{\pm\gamma}^m=W_{\lambda^m_{\pm\gamma}}=\Big(i\,
     \frac{d}{dt}\pm\gamma\Big)^m=\sum_{k=0}^m \binom{m}{k}i^k(\pm\gamma)^{m-k}\,\frac{d^k}{dt^k}\,.
\end{equation}
By the relations \eqref{eqn17} the mappings
\begin{equation*}
 \begin{aligned}
\mathbf{\Lambda}_{\gamma}^m &
:\widetilde{\mathbb{H}}^s_p(\mathbb{R}^+)
 \longrightarrow \widetilde{\mathbb{H}}^{s-m}_p(\mathbb{R}^+),
 \\[1ex]
\mathbf{\Lambda}_{-\gamma}^m
&:\mathbb{H}^s_p(\mathbb{R}^+)\longrightarrow\mathbb{H}^{s-m}_p(\mathbb{R}^+),
 \end{aligned}
\end{equation*}
are isomorphisms of the corresponding spaces if $\mathrm{Im}\,
\gamma>0$.

Taking into account formulae \eqref{eqn24} and \eqref{eqn25}, one
obtains the relation
\begin{align*}
\mathbf{\Lambda}^m_{\gamma}{\bf K}^1_c\varphi
&=\Big(i\,\frac{d}{dt}+\gamma\Big)^m{\bf K}^1_c\varphi=\sum_{k=0}^m
    \binom{m}{k}i^k\gamma^{m-k} \,\frac{d^k}{dt^k}\,{\bf K}^1_c\varphi\nonumber \\
& =\sum_{k=0}^m\binom{m}{k}i^k\gamma^{m-k}c^{-k}\Big({\bf K}^1_c\,
    \frac{d^k}{dt^k}\,\varphi\Big)(t)= \nonumber \\
& =c^{-m}{\bf K}^1_c\bigg(\sum_{k=0}^m
\binom{m}{k}i^k\left(c\,\gamma
    \right)^{m-k}\,\frac{d^k}{dt^k}\,\varphi\bigg)(t)= \nonumber \\
& =c^{-m}{\bf
K}^1_c\mathbf{\Lambda}^m_{c\,\gamma}\varphi,\qquad\varphi\in
    \widetilde{\mathbb{H}}^r_p(\mathbb{R}^+).
\end{align*}
Thus for $s=m=1,2,\ldots$, formula \eqref{eqn22} is proved.

If $s$ is a negative integer, $s=-1,-2,\ldots=-m$, formulae
\eqref{eqn22} can be established by applying the inverse operators
and $\mathbf{\Lambda}^{-m}_\gamma$ and
$\mathbf{\Lambda}^{-m}_{-c\gamma}$, respectively, from the left and
from the right to the already proven operator equality
 \[
\mathbf{\Lambda}^m_{\gamma}{\bf K}^1_c=c^{-m}{\bf
K}^1_c\mathbf{\Lambda}^m_{c\,\gamma},  \quad m=1,2,\ldots\,.
 \]
Thus one obtains
 \[
{\bf
K}^1_c\mathbf{\Lambda}^{-m}_{c\gamma}=c^{-m}\mathbf{\Lambda}^{-m}_\gamma{\bf
K}^1_c
     \quad {\rm or}\quad\mathbf{\Lambda}^{-m}_\gamma{\bf K}^1_c
     =c^m{\bf K}^1_c\mathbf{\Lambda}^{-m}_{c\gamma}
 \]
and for a negative $s=-1,-2,\ldots$, relation \eqref{eqn22} is also
proved.

In order to establish formula \eqref{eqn22} for non-integer values
of $s$, we can confine ourselves to the case $-1<s<0$. Indeed, any
non-integer value $s\in \mathbb{R}$ can be represented in the form
$s=s_0+m$, where $-1<s_0<0$ and $m$ is an integer. Therefore, if for
$s=s_0+m$ the operators in \eqref{eqn23} are correctly defined and
bounded, and if the relations in question are valid for $-1<s_0<0$,
then we can write
\begin{align*}
\mathbf{\Lambda}^s_{-\gamma}{\bf
K}^1_c=\mathbf{\Lambda}^{s_0+m}_{-\gamma}{\bf K}^1_c
     &=c^{-m}\mathbf{\Lambda}^{s_0}_{-\gamma}{\bf K}^1_c\mathbf{\Lambda}^{m}_{-c\,\gamma}
     =c^{-s_0-m}{\bf K}^1_c\mathbf{\Lambda}^{s_0}_{-c\,\gamma}
     \mathbf{\Lambda}^{m}_{-c\,\gamma}\nonumber\\
&=c^{-s_0-m}{\bf
K}^1_c\mathbf{\Lambda}^{s_0+m}_{-c\,\gamma}=c^{-s}{\bf K}^1_c
     \mathbf{\Lambda}^s_{-c\,\gamma}.
\end{align*}

Thus let us assume that $-1<s<0$ and consider the case
$0<\arg\,c<\pi$. Changing the orders of integration, we obtain
 \begin{equation}\label{eqn26}
 \begin{aligned}
\mathbf{\Lambda}^s_{-\gamma}{\bf
K}^1_c\varphi(t)&\,\,=\frac1{2\pi^2}\,r_+
    \int\limits_{-\infty}^\infty e^{-i\xi t}(\xi-\gamma)^s
    \int\limits_0^\infty e^{i\xi y}\int\limits_0^\infty
    \frac{\varphi(\tau)}{y-c\tau}\;d\tau\,dy\,d\xi\\[1ex]
&\,\,=\frac1{2\pi^2}\,r_+\int\limits_0^\infty\varphi(\tau)\int\limits_0^\infty
    \frac1{y-c\tau}\int\limits_{
    -\infty}^\infty e^{i\xi(y-t)}(\xi-\gamma)^sd\xi\,dy\,d\tau,
\end{aligned}
 \end{equation}
where $r_+$ is the restriction to $\mathbb{R}^+$. In order to study
the expression in the right-hand side of \eqref{eqn26}, one can use
a well known formula
\begin{equation*}
\begin{aligned}
\displaystyle
 &\int\limits_{-\infty}^\infty(\beta+ix)^{-\nu}e^{-ipx}\,dx=\begin{cases}
     0 \quad &\text{for} \quad p>0, \\
-\displaystyle\frac{2\pi(-p)^{\nu-1}e^{\beta\,p}}{\Gamma(\nu)}
     \quad&\text{for}\;\;p<0,\end{cases}\\[1.5ex]
& {\rm Re}\,\nu>0,\qquad{\rm Re}\beta>0 ,
\end{aligned}
\end{equation*}
 \cite[Formula
3.382.6]{GR94}. It can be rewritten in a more convenient form--viz.,
\begin{align}\label{eqn27}
&&\hskip-15mm\displaystyle\int\limits_{-\infty}^\infty
e^{i\mu\,\xi}(\xi
     -\gamma)^s\,d\xi=\left\{\begin{array}{ll}0 &\quad\text{ if } \;  \mu<0,\; {\rm Im}\,\gamma>0,\\
     \displaystyle\frac{2\pi\,\mu^{-s-1}e^{-\frac\pi2si+\mu\,\gamma i}}{
     \Gamma(-s)}&\quad\text{ if } \; \mu>0, \; {\rm Im}\,\gamma>0.\end{array}\right.
\end{align}

Applying \eqref{eqn27} to the last integral in \eqref{eqn26}, one
obtains
\begin{align}\label{eqn28}
\mathbf{\Lambda}^s_{-\gamma}{\bf K}^1_c\varphi(t)&\,\,
    =\displaystyle\frac{e^{-\frac\pi2si}}{\pi\Gamma(-s)}
    r_+\displaystyle\int\limits_0^\infty\varphi(\tau)\,d\tau\displaystyle
    \int\limits_t^\infty\displaystyle\frac{e^{i(y-t)\gamma}dy}{(y-t)^{1+s} (y-c\tau)}\nonumber\\
&=\displaystyle\frac{e^{-\frac\pi2si}}{\pi\Gamma(-s)}r_+\displaystyle
    \int\limits_0^\infty\varphi(\tau)\,d\tau\displaystyle\int\limits_0^\infty
    \displaystyle\frac{e^{i\gamma\,y}dy}{y^{1+s} (y+t-c\tau)}.
\end{align}

Let us also use the formula \cite[Formula 3.383.10]{GR94},
 \begin{equation}\label{eqn29}
  \begin{aligned}
\displaystyle
&\int\limits_0^\infty\frac{x^{\nu-1}e^{-\mu\,x}\,dx}{x+\beta}
     =\beta^{\nu-1}e^{\beta\,\mu}\Gamma(\nu)\Gamma(1-\nu,\beta\mu),\\[1ex]
&{\rm Re}\,\nu>0,\quad{\rm Re}\mu>0,\quad |\arg\,\beta|<\pi,
 \end{aligned}
 \end{equation}
and represent the operator \eqref{eqn28} in the form
\begin{equation}\label{eqn30}
\mathbf{\Lambda}^s_{-\gamma}{\bf K}^1_c\varphi(t)
    =\frac{e^{-\frac\pi2si}}\pi r_+\displaystyle\int\limits_0^\infty\displaystyle
    \frac{e^{-i\gamma(t-c\tau)}\Gamma(1+s,-i\gamma(t-c\tau))\varphi(\tau)
    \,d\tau}{(t-c\tau)^{1+s}}.
\end{equation}

Consider now the inverse composition ${\bf
K}^1_c\mathbf{\Lambda}^s_{-c\, \gamma} \varphi(t)$. Changing the
order of integration in the corresponding expression, one  obtains
\begin{align}\label{eqn31}
{\bf
K}^1_c\mathbf{\Lambda}^s_{-c\,\gamma}\varphi(t)&:=\frac1{2\pi^2}\,r_+
    \int\limits_0^\infty\frac1{t-c\,y}\int\limits_{-\infty}^\infty e^{-i\xi\,y}(\xi-c\,\gamma)^s
    \int\limits_0^\infty e^{i\xi\,\tau}\varphi(\tau)d\tau\,d\xi\,dy\nonumber\\
&=\frac1{2\pi^2}r_+\int\limits_0^\infty\varphi(\tau)\int\limits_0^\infty
    \frac1{t-c\,y}\int\limits_{-\infty}^\infty e^{i\xi(\tau-y)}
    (\xi-c\,\gamma)^sd\xi\,dy\,d\tau.
\end{align}
In order to compute the expression in the right-hand side of
\eqref{eqn31}, let us recall Formula~3.382.7 of \cite{GR94},
\begin{equation*}
\begin{aligned}
&\displaystyle\int\limits_{-\infty}^\infty(\beta-ix)^{-\nu}e^{-ipx}\,dx=\begin{cases}
     0\quad &\text{for}\quad p<0, \\
\displaystyle\frac{2\pi\,p^{\nu-1}e^{-\beta\,p}}{\Gamma(\nu)}
\quad&\text{for}
     \quad p>0,\end{cases}\\[1.5ex]
&{\rm Re}\,\nu>0,\quad{\rm Re}\beta>0,
\end{aligned}
\end{equation*}
and rewrite it in a form more suitable for our consideration--viz.,
\begin{equation}\label{eqn32}
\begin{aligned}
&\displaystyle\int\limits_{-\infty}^\infty
e^{i\mu\,\xi}(\xi+\omega)^s
     \,d\xi=\left\{\begin{array}{ll}0 &\quad\mu>0,\;{\rm Im}\,\omega>0,\\
     \displaystyle\frac{2\pi\,(-\mu)^{-s-1}e^{\frac\pi2si-\mu\,\omega i}}{
     \Gamma(-s)}&\quad\mu<0,\;{\rm
     Im}\,\omega>0,\end{array}\right.\\[1ex]
& {\rm Re}\,s<0,\quad\mu\in\mathbb{R},\quad\omega,\,s\in\mathbb{C}.
\end{aligned}
\end{equation}
Using \eqref{eqn32}, we represent \eqref{eqn31} in the form
 \begin{equation*}
 \begin{aligned}
{\bf K}^1_c\mathbf{\Lambda}^s_{-c\,\gamma}\varphi(t)
&=\displaystyle\frac{e^{\frac\pi2si}}{\pi\Gamma(-s)}\,r_+
    \displaystyle\int\limits_0^\infty \varphi(\tau)\,d\tau\displaystyle\int\limits_\tau^\infty
    \displaystyle\frac{e^{-ic\,\gamma(y-\tau)}\,dy}{(y-\tau)^{s+1}(t-c\,y)}\\
&=-\displaystyle\frac{e^{\frac\pi2si}}{\pi
c\Gamma(-s)}r_+\displaystyle
 \int\limits_0^\infty\varphi(\tau)\,d\tau\displaystyle
 \int\limits_0^\infty\displaystyle\frac{e^{-ic\gamma\,y}\,dy}{y^{s+1}(y-c^{-1}t+\tau)},
\end{aligned}
 \end{equation*}
and the application of formula \eqref{eqn29} leads to the
representation
\begin{align}\label{eqn33}
{\bf K}^1_c\mathbf{\Lambda}^s_{-c\,\gamma}\varphi(t)
     =&-\frac{c^{-1}e^{\frac\pi2si}}\pi
    r_+\int\limits_0^\infty\displaystyle\frac{e^{-ic\gamma(c^{-1}t-\tau)}\Gamma(1+s,
    -ic\gamma(c^{-1}t-\tau))\varphi(\tau)\,d\tau}{(\tau-c^{-1}t)^{1+s}}\nonumber\\
=& \frac{c^se^{-\frac\pi2si}}\pi r_+\int\limits_0^\infty
    \displaystyle\frac{e^{-i\gamma(t-c\,\tau)}\Gamma(1+s,-i\gamma(t-c\,\tau))
    \varphi(\tau)\,d\tau}{(t-c\,\tau)^{1+s}}.
\end{align}
Now the relations \eqref{eqn30} and  \eqref{eqn33} imply the
equality \eqref{eqn22} for $0<\arg\,c<\pi$.

In the case ${\rm Im}\,c=0$, $c<0$, we proceed as for
$0<\arg\,c<\pi$ and arrive at the formula
\begin{equation}\label{eqn34}
\mathbf{\Lambda}^s_{-\gamma}{\bf K}^1_{-1}\varphi(t)
    =\frac{e^{-\frac\pi2si}}\pi r_+\displaystyle\int\limits_0^\infty\displaystyle
    \frac{e^{-i\gamma(t+|c|\tau)}\Gamma(1+s,-i\gamma(t+|c|\tau))\varphi(\tau)
    \,d\tau}{(t+|c|\tau)^{1+s}},
\end{equation}
which is similar to \eqref{eqn30}. Further, instead of \eqref{eqn33}
we get
\begin{align}\label{eqn35}
{\bf K}^1_{-1}\mathbf{\Lambda}^s_{\gamma}\varphi(t)&=
     \frac{|c|^se^{\frac\pi2si}}\pi r_+\int\limits_0^\infty\displaystyle\frac{e^{-i\gamma
     (t+|c|\tau)}\Gamma(1+s,-i\gamma(t+|c|\tau))\varphi(\tau)\,d\tau}{(t+|c|\tau)^{1+s}},
\end{align}
and the relations \eqref{eqn34} and  \eqref{eqn35} lead to the
equality \eqref{eqn22} for ${\rm Im}\,c=0$, $c<0$. \hfill $\square$
 %
\begin{corollary}\label{c3.2}
Let $0<|\arg\,c|\leqslant\pi$, $\arg\,c\not=0$, $0<\arg \gamma<\pi$
and $-\pi<\arg(c\,\gamma)<0$. Then for arbitrary
$\gamma_0\in\mathbb{C}$ such that $0<\arg\,\gamma_0<\pi$ and
$-\pi<\arg(c\,\gamma_0)|<0$, one has
\begin{equation}\label{eqn36}
\mathbf{\Lambda}^s_{-\gamma}{\bf
K}^1_c=c^{-s}W_{g_{-\gamma,-\gamma_0}}{\bf K}^1_c
   \mathbf{\Lambda}^s_{-c\,\gamma_0},
\end{equation}
where
\begin{align}\label{eqn37}
g^s_{-\gamma,-\gamma_0}(\xi):=\left(\frac{\xi-\gamma}{\xi-\gamma_0}\right)^s.
\end{align}

If, in addition, $1<p<\infty$ and $1/p-1<r<1/p$ then equality
\eqref{eqn36} can be supplemented as follows
\begin{align}\label{eqn38}
\mathbf{\Lambda}^s_{-\gamma}{\bf K}^1_c=c^{-s}\left[{\bf
K}^1_cW_{g^s_{-\gamma,-\gamma_0}}
   +\mathbf{T}\right]\mathbf{\Lambda}^s_{-c\,\gamma_0},
\end{align}
where
$\mathbf{T}\;:\;\widetilde{\mathbb{H}}{}^r_p(\mathbb{R}^+)\to\mathbb{H}^r_p(\mathbb{R}^+)$
is a compact operator, and if $c$ is a real negative number, then
$c^{-s}:=|c|^{-s}e^{-\pi si}$.
\end{corollary}
\textbf{Proof.}  It follows from equalities \eqref{eqn16} and
\eqref{eqn22} that
\begin{eqnarray*}
\mathbf{\Lambda}^s_{-\gamma}{\bf K}^1_c=\mathbf{\Lambda}^s_{-\gamma}
     \mathbf{\Lambda}^{-s}_{-\gamma_0}\mathbf{\Lambda}^s_{-\gamma_0}{\bf K}^1_c
     =c^{-s}W_{g_{-\gamma,-\gamma_0}}{\bf K}^1_c\mathbf{\Lambda}^s_{-c\,\gamma_0}
\end{eqnarray*}
and \eqref{eqn36} is proved. If $1<p<\infty$ and $1/p-1<r<1/p$, then
the commutator
 \[
\mathbf{T}:=W_{g^s_{-\gamma,-\gamma_0}}{\bf K}^1_c-{\bf
K}^1_cW_{g^s_{-\gamma,-\gamma_0}}
     \;:\;\widetilde{\mathbb{H}}{}^r_p(\mathbb{R}^+)\to\mathbb{H}^r_p(\mathbb{R}^+)
 \]
of Mellin and Fourier convolution operators is correctly defined and
bounded. It is compact for $r=0$ and all $1<p<\infty$ (see
\cite{Du74,Du87}). Due to Krasnoselsky interpolation theorem (see
\cite{Kr60} and also \cite[Sections 1.10.1 and 1.17.4]{Tr95}), the
operator $\mathbf{T}$ is compact in all $\mathbb{L}_r$-spaces for
$1/p-1<r<1/p$. Therefore, the equality \eqref{eqn36}, can be
rewritten as
\begin{eqnarray*}
\mathbf{\Lambda}^s_{-\gamma}{\bf K}^1_c=c^{-s}\left[{\bf
K}^1_cW_{g^s_{-\gamma,-\gamma_0}}
     +\mathbf{T}\right]\mathbf{\Lambda}^s_{-c\,\gamma_0}\, ,
\end{eqnarray*}
and we are done \rbx
 %
\begin{remark}\label{r1.7}
The assumption $1/p-1<r<1/p$ in \eqref{eqn38} cannot be relaxed.
Indeed, the operator $W_{g^s_{-\gamma,-\gamma_0}}{\bf
K}^1_c=\mathbf{\Lambda}^s_{-\gamma}
\mathbf{\Lambda}^{-s}_{-\gamma_0}{\bf K}^1_c\;:\;\widetilde{
\mathbb{H}}{}^r_p(\mathbb{R}^+)\to\mathbb{H}^r_p(\mathbb{R}^+)$ is
bounded for all $r\in\mathbb{R}$ (see \eqref{eqn23}). But the
operator ${\bf K}^1_cW_{g^s_{-\gamma,
-\gamma_0}}\;:\;\widetilde{\mathbb{H}}{}^r_p
(\mathbb{R}^+)\to\mathbb{H}^r_p(\mathbb{R}^+)$ is bounded only for
$1/p-1<r<1/p$ because the function $g^s_{-\gamma,-\gamma_0}(\xi)$
has an analytic extension into the lower half-plane but not into the
upper one.
\end{remark}

 %
\section{Algebra Generated by Mellin and Fourier Convolution Operators}
\label{s4}

In the present section we recall some results on Banach algebra,
generated by Fourier and Mellin convolution operators in the
Lebesgue space with weight from \cite{Du87}, revised in \cite{Du13}.
The exposition follows \cite[Section 2]{Du13}. For more general
algebras we refer the reader to \cite{Du87} and to \cite{Du74,Th85}.

Let us consider the Banach algebra $\mathfrak{A}_p(\mathbb{R}^+)$
generated by Mellin convolution and Fourier convolution operators in
the Lebesgue space $\mathbb{L}_p(\mathbb{R}^+)$. In particular, this
algebra contains the operators
\begin{equation}\label{eqn39}
    \mathbf{A}:=\sum_{j=1}^m\mathfrak{M}^0_{a_j}W_{b_j},
\end{equation}
and their compositions. Here $\mathfrak{M}^0_{a_j}$ are Mellin
convolution operators with continuous $N\times N$ matrix symbols
$a_j\in C\mathfrak{M}_p(\overline{R})$, $W_{b_j}$ are Fourier
convolution operators with $N\times N$ matrix symbols $b_j\in
C\mathfrak{M}_p(\overline{\mathbb{R}}\setminus\{0\}):=C\mathfrak{M}_p(\overline{\mathbb{R}}^-
\cup\overline{\mathbb{R}}^+)$. The algebra of $N\times N$ matrix
$\mathbb{L}_p$-multipliers
$C\mathfrak{M}_p(\overline{\mathbb{R}}\setminus\{0\})$ consists of
those piecewise-continuous $N\times N$ matrix multipliers
$b\in\mathfrak{M}_p(\mathbb{R})\cap PC(\overline{\mathbb{R}})$ which
are continuous on the semi-axes $\mathbb{R}^-$ and $\mathbb{R}^+$
but might have finite jump discontinuities at $0$ and at the
infinity.

Note that the algebra $\mathfrak{A}_p(\mathbb{R}^+)$ is actually a
subalgebra of the Banach algebra $\mathfrak{F}_p(\mathbb{R}^+)$
generated by the Fourier convolution operators $W_a$ with
piecewise-constant symbols $a(\xi)$ in the space
$\mathbb{L}_p(\mathbb{R}^+)$. Let
$\mathfrak{S}(\mathbb{L}_p(\mathbb{R}^+))$ denote the ideal of all
compact operators in $\mathbb{L}_p(\mathbb{R}^+)$. Since in the
scalar case $N=1$ the quotient algebra
$\mathfrak{F}_p(\mathbb{R}^+)/ \mathfrak{S}
(\mathbb{L}_p(\mathbb{R}^+))$ is commutative, the following
proposition is true.

\begin{figure}[ht]
\centering
\includegraphics[height=49mm]{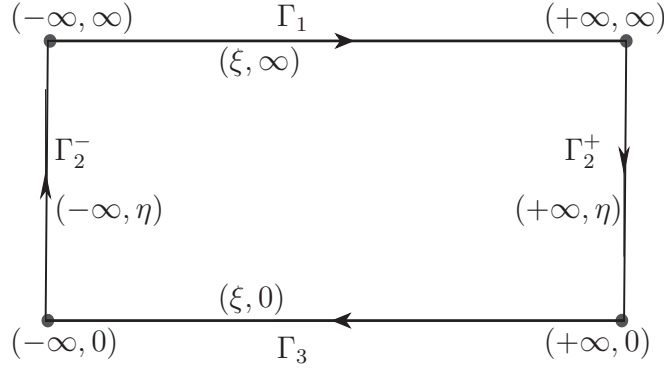}
\caption{The domain $\mathfrak{R}$ of definition of the symbol
$\mathcal{A}_p(\xi, \eta)$.} \label{Fig:1}
\end{figure}
 %
\begin{proposition}[\cite{Du87} and {\cite[Corollary 3.10]{Du13}}]\label{p4.1}
If $N=1$, then the quotient algebra $\mathfrak{A}_p(\mathbb{R}^+)/
\mathfrak{S}(\mathbb{L}_p (\mathbb{R}^+))$ is commutative.
\end{proposition}

To describe the symbol of the operator $\boldsymbol{A}$ of
\eqref{eqn39}, consider the infinite clockwise oriented
``rectangle'' $\mathfrak{R}:=\Gamma_1\cup\Gamma_2^-\cup\Gamma_2^+
\cup\Gamma_3$, where (cf. Figure~1)
 \[
 \Gamma_1:=\overline{\mathbb{R}}\times\{+\infty\},\;\;\Gamma^\pm_2:=\{\pm\infty\}
     \times\overline{\mathbb{R}}^+,\;\;\Gamma_3:=\overline{\mathbb{R}}\times\{0\}.
 \]

\noindent The symbol $\mathcal{A}_p(\omega)$ of the operator
$\boldsymbol{A}$ in \eqref{eqn39} is a function on the set
$\mathfrak{R}$, viz.
\begin{equation}\label{eqn40}
    \mathcal{A}_p(\omega):=\begin{cases}
        \displaystyle \sum_{j=1}^m a_j(\xi)(b_j)_p(\infty,\xi), & \omega=(\xi,\infty)\in\overline{\Gamma}_1, \\
        \displaystyle \sum_{j=1}^m a_j(+\infty)b_j(-\eta), & \omega=(+\infty,\eta)\in\Gamma^+_2, \\
        \displaystyle \sum_{j=1}^m a_j(-\infty)b_j(\eta), & \omega=(-\infty,\eta)\in\Gamma^-_2, \\
        \displaystyle \sum_{j=1}^m (a_j)_p(\infty,\xi)(b_j)_p(0,\xi), & \omega=(\xi,0)\in\overline{\Gamma}_3.
                \end{cases}
\end{equation}
In \eqref{eqn40} for a piecewise continuous function $g\in
PC(\overline{\mathbb{R}})$ we use the notation
\begin{equation}\label{eqn41}
\begin{aligned}
    g_p(\infty,\xi) & :=\frac12\,\big[g(+\infty)+g(-\infty)\big]- \\
    &\qquad\qquad -\frac i2\,\big[g(+\infty)-g(-\infty)\big]
            \cot\pi\Big(\frac1p-i\xi\Big), \\
    g_p(t,\xi) & :=\frac12\,\big[g(t+0)+g(t-0)\big]- \\
    &\qquad\qquad -\frac i2\,\big[g(t+0)-g(t-0)\big]\cot\pi\Big(\frac1p-i\xi\Big),
\end{aligned}
\end{equation}
where $t,\xi\in\mathbb{R}$.

\begin{figure}[ht]
\centering
\includegraphics[height=25mm]{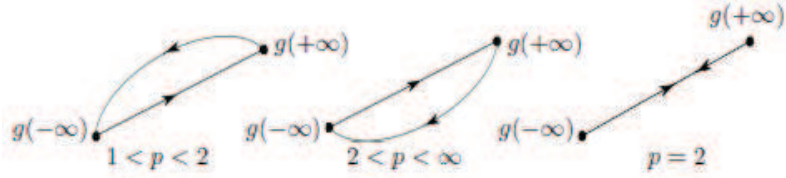}
\caption{Arc condition.} \label{Fig:2}
\end{figure}

\noindent {\bf Arc condition} (\cite{GK79,Wi60}): The function
$g_p(\infty,\xi)$  connects the point $g(-\infty)$ with
$g(+\infty)$. More precisely, it fills up the discontinuity of the
function $g$ at $\infty$ with an oriented arc of the circle such
that from every point of the arc the oriented interval
$[g(-\infty),g(+\infty)]$ is seen under the angle $\pi/p$. Moreover,
the oriented arc lies on the left of the oriented interval if
$1/2<1/p<1$ (i.e., if $1<p<2$) and the oriented arc is on the right
of the oriented interval if $0<1/p<1/2$ (i.e., if $2<p<\infty$). For
$p=2$ the oriented arc coincides with the oriented interval (see
Figure 2).

A similar geometric interpretation is valid for the function
$g_p(t,\xi)$, which connects the points $g(t-0)$ and $g(t+0)$ when
$g$ has a jump at $t\in\overline{\mathbb{R}}$.

The image of the function $\det\mathcal{A}_p(\omega)$, $\omega\in
\mathfrak{R}$ is a closed curve in the complex plane. It follows
from the continuity of the symbol at the angular points of the
rectangle $\mathfrak{R}$ where the one-sided limits coincide. Thus
\begin{align*}
    \mathcal{A}_p(\pm\infty,\infty) & =\sum_{j=1}^m a_j(\pm\infty)b_j(\mp\infty), \\
    \mathcal{A}_p(\pm\infty,0) & =\sum_{j=1}^m a_j(\pm\infty)b_j(0\mp0).
\end{align*}
Hence, if the symbol of the corresponding operator is elliptic, i.e.
if
 \[
    \inf_{\omega\in\mathfrak{R}} \big|\det\mathcal{A}_p(\omega)\big|>0,
 \]
the increment of the argument $(1/2\pi)\arg \mathcal{A}_p(\omega)$
when $\omega$ ranges through $\mathfrak{R}$ in the direction of
orientation, is an integer. It is called the winding number or the
index of the curve $\Gamma:=\{z\in \mathbb{C}: z=
\det\mathcal{A}_p(\omega),\;\omega \in \mathfrak{R}\}$ and is
denoted by ${\rm ind}\,\det\mathcal{A}_p$.
 %
\begin{theorem}[{\cite[Theorem 3.13]{Du13}}]\label{t4.2}
Let $1<p<\infty$ and let ${\bf A}$ be defined by \eqref{eqn39}. The
operator ${\bf
A}\;:\;\mathbb{L}_p(\mathbb{R}^+)\longrightarrow\mathbb{L}_p(\mathbb{R}^+)$
is Fredholm if and only if its symbol $\mathcal{A}_p(\omega)$ is
elliptic. If $\mathbf{A}$ is Fredholm, then the index  of this
operator is
 \[
    {\rm Ind}\mathbf{A}=-{\rm ind}\det\mathcal{A}_p.
 \]
\end{theorem}

If $\mathcal{A}_p(\omega)$ is the symbol of an operator $\mathbf{A}$
in \eqref{eqn39}, then the set
$\mathcal{R}(\mathcal{A}_p):=\{\mathcal{A}_p(\omega)\in\mathbb{C}:\omega\in
\mathfrak{R}\}$ coincides with the essential spectrum of
$\mathbf{A}$. Recall that the essential spectrum
$\sigma_{ess}(\mathbf{A})$ of a bounded operator $\mathbf{A}$ is the
set of all $\lambda\in \mathbb{C}$ such that the operator
$\mathbf{A}-\lambda I$ is not Fredholm in
$\mathbb{L}_p(\mathbb{R}^+)$ or, equivalently, the coset
$[\mathbf{A}-\lambda I]$ is not invertible in the quotient algebra
$\mathfrak{A}_p(\mathbb{R}^+)/
\mathfrak{S}(\mathbb{L}_p(\mathbb{R}^+))$. Then, due to Banach
theorem, the essential norm $\|\!|\mathbf{A}\|\!|$ of the operator
$\mathbf{A}$ can be estimated as follows
\begin{equation}\label{eqn42}
\sup_{\omega\in\omega}|\mathcal{A}_p(\omega)|\leqslant\|\!|\mathbf{A}\|\!|
     :=\inf_{\mathbf{T}\in\mathfrak{S}(\mathbb{L}_p(\mathbb{R}^+))}
     \big\|(\mathbf{A}+\mathbf{T})\mid\mathcal{L}(\mathbb{L}_p(\mathbb{R}^+))\big\|.
\end{equation}
The inequality \eqref{eqn42} enables one to extend continuously the
symbol map \eqref{eqn40}
 \[
[\mathbf{A}]\longrightarrow\mathcal{A}_p(\omega),\quad
[\mathbf{A}]\in\mathfrak{A}_p(
     \mathbb{R}^+)/\mathfrak{S}(\mathbb{L}_p(\mathbb{R}^+))
 \]
on the whole Banach algebra $\mathfrak{A}_p(\mathbb{R}^+)$. Now,
applying Theorem \ref{t4.2} and a standard methods, cf.
\cite[Theorem 3.2]{Du87}, one can derive the following result.
 %
\begin{corollary}[{\cite[Corollary 3.15]{Du13}}]\label{c4.3}
Let $1<p<\infty$ and $\mathbf{A}\in\mathfrak{A}_p(\mathbb{R}^+)$.
The operator
$\mathbf{A}:\mathbb{L}_p(\mathbb{R}^+)\longrightarrow\mathbb{L}_p(\mathbb{R}^+)$
is Fredholm if and only if its symbol $\mathcal{A}_p(\omega)$ is
elliptic. If $\mathbf{A}$ is Fredholm, then
 $$
{\rm Ind}\,\mathbf{A}=-{\rm ind}\,\mathcal{A}_p.
 $$
\end{corollary}

 %
\section{Fredholm properties of Mellin Convolution Operators in the Bessel Potential Spaces.}
\label{s5}

As it was already mentioned, the primary aim of the present paper is
to study Fredholm properties and the invertibility of Mellin
convolution operators $\mathfrak{M}^0_a$ acting in Bessel potential
spaces, namely,
\begin{equation*}
\mathfrak{M}^0_a:\widetilde{\mathbb{H}}^s_p(\mathbb{R}^+)\longrightarrow
     \mathbb{H}^s_p(\mathbb{R}^+).
\end{equation*}
The symbols of these operators are $N\times N$ matrix functions
$a\in C\mathfrak{M}^0_p (\overline{\mathbb{R}}),$ continuous on the
real axis $\mathbb{R}$ with the only possible jump at infinity.
 %
\begin{theorem}\label{t5.1}
Let $s\in\mathbb{R}$ and $1<p<\infty$.
\begin{enumerate}
\item If the conditions of Theorem \ref{t3.1} hold, then the Mellin
convolution operator $\mathbf{K}^1_c$,
\begin{align}\label{eqn43}
\mathbf{K}^1_c\;:\;\widetilde{\mathbb{H}}{}^r_p(\mathbb{R}^+)\to\mathbb{H}^r_p(\mathbb{R}^+)
\end{align}
is lifted to the equivalent operator
\begin{align*}
\mathbf{\Lambda}^s_{-\gamma}\mathbf{K}^1_c\mathbf{\Lambda}^{-s}_{\gamma}
     =c^{-s}\mathbf{K}^1_cW_{g^s_{-c\gamma,\gamma}}
     :\mathbb{L}_p(\mathbb{R}^+)\to\mathbb{L}_p(\mathbb{R}^+),
\end{align*}
where $c^{-s}=|c|^{-s}e^{-is\arg\,c}$ and the function
$g^s_{-c\gamma,\gamma}$ is defined in \eqref{eqn37}.

\item If conditions of Corollary \ref{c3.2} hold, the Mellin
convolution operator between Bessel potential spaces \eqref{eqn43}
is lifted to the equivalent operator
\begin{align*}
\mathbf{\Lambda}^s_{-\gamma}\mathbf{K}^1_c\mathbf{\Lambda}^{-s}_{\gamma}
     &=\!c^{-s}W_{g^s_{-\gamma,-\gamma_0}}\mathbf{K}^1_cW_{g^s_{-c\gamma_0,\gamma}}\\
     & =\!c^{-s}
     \mathbf{K}^1_cW_{g^s_{-\gamma,-\gamma_0}g^s_{-c\gamma_0,\gamma}}\!+\!\mathbf{T}
:\mathbb{L}_p(\mathbb{R}^+)\to\mathbb{L}_p(\mathbb{R}^+),\nonumber
\end{align*}
where $\mathbf{T}\;:\;\mathbb{L}_p(\mathbb{R}^+)
\to\mathbb{L}_p(\mathbb{R}^+)$ is a compact operator.
 \end{enumerate}
\end{theorem}
\textbf{Proof.}  By Theorem \ref{t1.1}, using the lifting procedure,
one obtains the following equivalent operator
 \[
\mathbf{\Lambda}^s_{-\gamma}\mathbf{K}^1_c\mathbf{\Lambda}^{-s}_{\gamma}
     \;:\;\mathbb{L}_p(\mathbb{R}^+)\to\mathbb{L}_p(\mathbb{R}^+).
 \]
In order to proceed, we need two formulae
 \begin{align}\label{eqn44}
\mathbf{\Lambda}^s_{-c\gamma}\mathbf{\Lambda}^{-s}_{\gamma}=W_{g^s_{-c\gamma,\gamma}},
     \qquad W_{g^s_{-\gamma,-\gamma_0}}W_{g^s_{-c\gamma_0,\gamma}}
     =W_{g^s_{-\gamma,-\gamma_0}g^s_{-c\gamma_0,\gamma}}.
 \end{align}
The first relation holds because, by the conditions of Theorem
\ref{t3.1}, $0<\arg\,\gamma<\pi$ and the second one holds because
$g^s_{-\gamma,-\gamma_0}(\xi)$ has a smooth, uniformly bounded
analytic extension in the complex lower half plane.

If conditions of Theorem \ref{t3.1} are satisfied, we use the
relations \eqref{eqn22}, \eqref{eqn44}. Thus
 \[
\mathbf{\Lambda}^s_{-\gamma}\mathbf{K}^1_c\mathbf{\Lambda}^{-s}_{\gamma}
     =c^{-s}\mathbf{K}^1_c\mathbf{\Lambda}^s_{-c\gamma}\mathbf{\Lambda}^{-s}_{\gamma}
     =c^{-s}\mathbf{K}^1_cW_{g^s_{-c\gamma,\gamma}}.
 \]
If conditions of Corollary \ref{c3.2} hold, we successively apply
formulae \eqref{eqn36}, \eqref{eqn38}, both formulae \eqref{eqn44},
so that
\begin{eqnarray*}
\mathbf{\Lambda}^s_{-\gamma}\mathbf{K}^1_c\mathbf{\Lambda}^{-s}_{\gamma}
     &=&c^{-s}W_{g^s_{-\gamma,-\gamma_0}}\mathbf{K}^1_c\mathbf{\Lambda}^s_{-c\gamma}
     \mathbf{\Lambda}^{-s}_{\gamma}\nonumber\\[2mm]
&=&c^{-s}W_{g^s_{-\gamma,-\gamma_0}}\mathbf{K}^1_c
     W_{g^s_{-c\gamma_0,\gamma}} =c^{-s}\mathbf{K}^1_cW_{g^s_{-\gamma,-\gamma_0}}W_{g^s_{-c\gamma_0,\gamma}}+\mathbf{T},
\end{eqnarray*}
which completes the proof. \rbx \vskip5mm
 %
 \begin{remark}
The operator $\mathbf{K}^1_1$ is the Hilbert transform
$\mathbf{K}^1_1=-\pi i S_{\mathbb{R}^+}=\pi iW_{\rm sign}$ and does
not satisfy the condition $\arg\,c\not=0$ of Theorem \ref{t5.1}. As
already emphasized in Theorem~\ref{t1.1}, this case is essentially
different. Considered as acting between the Bessel potential spaces
\eqref{eqn43}, $\mathbf{K}^1_1$ is lifted to the equivalent Fourier
convolution operator
\begin{align*}
\mathbf{\Lambda}^s_{-\gamma}\mathbf{K}^1_1\mathbf{\Lambda}^{-s}_{\gamma}
     =W_{\pi ig^s_{-\gamma,\gamma}{\rm sign}}
     \;:\;\mathbb{L}_p(\mathbb{R}^+)\to\mathbb{L}_p(\mathbb{R}^+),
\end{align*}
cf. Theorem \ref{t1.1}.
 \end{remark}
 %
\begin{theorem}\label{t5.2}
Let $c_j, d_j\in\mathbb{C}$, $-\pi\leqslant\arg\,c_j<\pi$
$\arg\,c_j\not=0$, for $j=1,\ldots,n$, $0<\arg \gamma<\pi$,
$-\pi<\arg(c_j\gamma)<0$ for $j=1,\ldots,m$  and
$0<\arg(c_j\gamma)<\pi$ for $j=m+1,\ldots,n$. The Mellin convolution
operator $\mathbf{A}$,
\begin{align*}
{\bf
A}=\sum_{j=1}^nd_j\mathbf{K}^1_{c_j}\;:\;\widetilde{\mathbb{H}}{}^r_p(
     \mathbb{R}^+)\to\mathbb{H}^r_p(\mathbb{R}^+),
\end{align*}
is lifted to the equivalent operator
\begin{align}\label{e69a}
\hskip-15mm\mathbf{\Lambda}^s_{-\gamma}{\bf
A}\mathbf{\Lambda}^{-s}_{\gamma}
    & =\sum_{j=0}^md_jc^{-s}_j\mathbf{K}^1_{c_j}W_{g^s_{-c_j\gamma,-\gamma}}
     +\sum_{j=m+1}^nd_jc^{-s}_jW_{g^s_{-\gamma,-\gamma_j}}
     \mathbf{K}^1_{c_j}W_{g^s_{-c_j\gamma_j,\gamma}}\\
&=\sum_{j=0}^md_jc^{-s}_j\mathbf{K}^1_{c_j}W_{g^s_{-c_j\gamma,\gamma}}
     +\sum_{j=m+1}^nd_jc^{-s}_j\mathbf{K}^1_{c_j}W_{g^s_{-\gamma,-\gamma_j}
     g^s_{-c_j\gamma_j,\gamma}}+\mathbf{T}\nonumber
\end{align}
in the $\mathbb{L}_p(\mathbb{R}^+)$ space, where
$c^{-s}=|c|^{-s}e^{-is\arg\,c}$ and $\gamma_j$ are such that $0<\arg
\gamma_j<\pi$, $-\pi<\arg(c_j\,\gamma_j)<0$ for $j=m+1,\ldots,n$.
$\mathbf{T}\;:\;\mathbb{L}_p(\mathbb{R}^+)
\to\mathbb{L}_p(\mathbb{R}^+)$ is a compact operator.
\end{theorem}
\textbf{Proof.}  The proof is a direct consequence of Theorem
\ref{t5.1}. \rbx
 %
\begin{theorem}\label{t5.3}
Let $s\in\mathbb{R}$ and $1<p<\infty$. If conditions of Theorem
\ref{t3.1} hold, then the Mellin convolution operator
$\mathbf{K}^2_c$,
\begin{align}\label{eqn45}
\mathbf{K}^2_c\;:\;\widetilde{\mathbb{H}}{}^r_p(\mathbb{R}^+)\to\mathbb{H}^r_p(\mathbb{R}^+)
\end{align}
is lifted to the equivalent operator
\begin{align}\label{eqn46}
\mathbf{\Lambda}^s_{-\gamma}\mathbf{K}^2_c\mathbf{\Lambda}^{-s}_{\gamma}
     =c^{-s}\left[\mathbf{K}^2_c-sc^{-1}\mathbf{K}^1_c\right]W_{g^s_{-c\gamma,\gamma}}
     +s\,\gamma\,c^{-s}\mathbf{K}^1_cW_{g^{s-1}_{-c\gamma,\gamma}}\mathbf{\Lambda}^{-1}_{\gamma}
\end{align}
in $\mathbb{L}_p(\mathbb{R}^+)$ space, where
$c^{-s}=|c|^{-s}e^{-is\arg\,c}$, the function
$g^s_{-c\gamma,\gamma}$ is defined in \eqref{eqn37}, and the last
summand in \eqref{eqn46}, namely, the operator
 \begin{align}\label{eqn47}
\mathbf{T}:=s\,\gamma\,c^{-s}\mathbf{K}^1_cW_{g^{s-1}_{-c\gamma,\gamma}}\mathbf{\Lambda}^{-1}_{\gamma}
     \;:\;\mathbb{L}_p(\mathbb{R}^+)\to\mathbb{L}_p(\mathbb{R}^+),
 \end{align}
is compact. Moreover, if conditions of Corollary \ref{c3.2} hold,
the Mellin convolution operator $\mathbf{K}^2_c$ between Bessel
potential spaces \eqref{eqn45} is lifted to the equivalent operator
\begin{align}\label{eqn48}
\mathbf{\Lambda}^s_{-\gamma}\mathbf{K}^2_c\mathbf{\Lambda}^{-s}_{\gamma}
     =& c^{-s}W_{g^s_{-\gamma,-\gamma_0}}\left[\mathbf{K}^2_c-sc^{-1}\mathbf{K}^1_c\right]
     W_{g^s_{-c\gamma_0,\gamma}}\nonumber\\[1ex]
&+s\,\gamma\,c^{-s}W_{g^s_{-\gamma,-\gamma_0}}\mathbf{K}^1_cW_{g^{s-1}_{-c\gamma_0,\gamma}}
     \mathbf{\Lambda}^{-1}_{\gamma}\nonumber\\[1ex]
=&c^{-s}\left[\mathbf{K}^2_c-sc^{-1}\mathbf{K}^1_c\right]W_{g^s_{-\gamma,-\gamma_0}
    g^s_{-c\,\gamma_0,\gamma}}+\mathbf{T}_0
\end{align}
in $\mathbb{L}_p(\mathbb{R}^+)$ space, and the operator
$\mathbf{T}_0\;:\;\mathbb{L}_p(\mathbb{R}^+)\to
\mathbb{L}_p(\mathbb{R}^+)$ is compact.
\end{theorem}
\textbf{Proof.}  If the conditions of Theorem \ref{t3.1} are
satisfied, then  ${\rm Im}\,\gamma>0$ and ${\rm Im}\,c\,\gamma<0$.
Hence
 \[
\frac1{(t-c)^2}=\lim_{\varepsilon\to0}\frac1{2\varepsilon
i}\left[\frac1{t-c-\varepsilon i}
     -\frac1{t-c+\varepsilon i}\right]
 \]
and we have
\begin{align*}
\mathbf{\Lambda}^s_{-\gamma}\mathbf{K}^2_c\mathbf{\Lambda}^{-s}_{\gamma}&=
\lim_{\varepsilon\to0}\frac1{2\varepsilon
i}\mathbf{\Lambda}^s_{-\gamma} \left[\mathbf{K}^1_{c +\varepsilon
i}-\mathbf{K}^1_{c-\varepsilon i}\right]
\mathbf{\Lambda}^{-s}_{\gamma}\displaybreak[2]\\[1ex]
&=\lim_{\varepsilon\to0}\frac1{2\varepsilon i} \left[(c+\varepsilon
i)^{-s}\mathbf{K}^1_{c +\varepsilon
i}\mathbf{\Lambda}^s_{-(c+\varepsilon i)\gamma} - (c-\varepsilon
i)^{-s}\mathbf{K}^1_{c-\varepsilon
i}\mathbf{\Lambda}^s_{-(c-\varepsilon
i)\gamma}\right]\mathbf{\Lambda}^{-s}_{\gamma}\displaybreak[2]\\[1ex]
&=\lim_{\varepsilon\to0}\left\{\frac{(c+\varepsilon
i)^{-s}-(c-\varepsilon i)^{-s}}{2\varepsilon i}\mathbf{K}^1_{c
+\varepsilon i}\mathbf{\Lambda}^s_{-(c+\varepsilon
i)\gamma}\right.\displaybreak[2]\\[1ex]
  &\quad \left.- (c-\varepsilon
i)^{-s}\frac1{2\varepsilon i}\left[ \mathbf{K}^1_{c+\varepsilon
i}-\mathbf{K}^1_{c-\varepsilon i}\right]\mathbf{\Lambda}^s_{-(c-\varepsilon i)\gamma}\right.\\[1ex]
&\quad \left.- (c-\varepsilon i)^{-s}\mathbf{K}^1_{c-\varepsilon
i}\frac1{2\varepsilon i}\left[ \mathbf{\Lambda}^s_{-(c+\varepsilon
i)\gamma}-\mathbf{\Lambda}^s_{-(c-\varepsilon
i)\gamma}\right]\right\}\mathbf{\Lambda}^{-s}_{\gamma}\\[1ex]
&=-s\,c^{-s-1}\mathbf{K}^1_c\mathbf{\Lambda}^s_{-c\,\gamma}\mathbf{\Lambda}^{-s}_{\gamma}
+c^{-s}\mathbf{K}^2_c\mathbf{\Lambda}^s_{-c\,\gamma}\mathbf{\Lambda}^{-s}_{\gamma}\\[1ex]
&\quad + c^{-s}\mathbf{K}^1_c\lim_{\varepsilon\to0}
\mathcal{F}^{-1}\frac{(\xi-c\,\gamma-\varepsilon\gamma
i)^s-(\xi-c\,\gamma+\varepsilon\gamma i)^s}{2\varepsilon
i}\mathcal{F}\mathbf{\Lambda}^{-s}_{\gamma}\\[1ex]
&=c^{-s}\left[\mathbf{K}^2_c-sc^{-1}\mathbf{K}^1_c\right]W_{g^s_{-c\gamma,\gamma}}
+s\,\gamma\,c^{-s}\mathbf{K}^1_c\mathbf{\Lambda}^{s-1}_{-c\,\gamma}\mathbf{\Lambda}^{-s}_{\gamma}\\[1ex]
&=c^{-s}\left[\mathbf{K}^2_c-sc^{-1}\mathbf{K}^1_c\right]W_{g^s_{-c\gamma,\gamma}}
+s\,\gamma\,c^{-s}\mathbf{K}^1_cW_{g^{s-1}_{-c\gamma,\gamma}}\mathbf{\Lambda}^{-1}_{\gamma}.
 \end{align*}
Thus formula \eqref{eqn46} is proved. To verify the compactness of
the operator $\mathbf{T}$ in \eqref{eqn47}, let us rewrite it as
follows
\begin{align}\label{eqn49}
\mathbf{T}&=s\,\gamma\,c^{-s}\mathbf{K}^1_cW_{g^{s-1}_{-c\gamma,\gamma}}\mathbf{\Lambda}^{-1}_{\gamma}
\nonumber\\[1ex]
   & =s\,\gamma\,c^{-s}(1-h)\mathbf{K}^1_cW_{g^{s-1}_{-c\gamma,\gamma}}\mathbf{\Lambda}^{-1}_{\gamma}
    +s\,\gamma\,c^{-s}h\mathbf{K}^1_cW_{g^{s-1}_{-c\gamma,\gamma}}\mathbf{\Lambda}^{-1}_{\gamma},
\end{align}
where $h\in C^\infty(\mathbb{R})^+$ is a smooth function having a
compact support and equal to $1$ in a neighborhood of $0$. Since
$1-h(t)$ vanishes in the neighbourhood of $0$, the operator
$(1-h)\mathbf{K}^1_c$ has a smooth kernel and is compact in
$\mathbb{L}_p(\mathbb{R}^+)$. The second summand in \eqref{eqn49} is
compact since $h$ commutes with the Mellin $\mathbf{K}^1_c$ and
Fourier $W_{g^{s-1}_{-c\gamma,\gamma}}$ convolutions modulo compact
operators, i.e.
 \[
s\,\gamma\,c^{-s}h\mathbf{K}^1_cW_{g^{s-1}_{-c\gamma,\gamma}}\mathbf{\Lambda}^{-1}_{\gamma}
      =s\,\gamma\,c^{-s}\mathbf{K}^1_cW_{g^{s-1}_{-c\gamma,\gamma}}h\mathbf{\Lambda}^{-1}_{\gamma}+\mathbf{T}_1,
 \]
where $\mathbf{T}_1$ is compact in $\mathbb{L}_p(\mathbb{R}^+)$ (see
Proposition \ref{p4.1} and \cite[Lemma 7.4]{Du79}, \cite[Lemma
1.2]{Du87}). Note, that due to the Sobolev embedding theorem,  the
operator $h\mathbf{\Lambda}^{-1}_{\gamma}$ is also compact in
$\mathbb{L}_p(\mathbb{R}^+)$, because ${\rm supp}\,h$ is compact.
Finally, formula \eqref{eqn48} can be derived from \eqref{eqn46}
similarly to considerations of Theorem \ref{t5.1}. \rbx
 \sloppy
\begin{remark}\label{r5.3a}
The operators $\mathbf{K}^n_c$, $n=3,4,\ldots$, can be treated
analogously to the approach of Corollary \ref{t5.3}. Indeed, let us
represent the operator $\mathbf{K}^n_c$ in the form
 $$
\mathbf{K}^n_c\varphi=\lim_{\varepsilon\to0}\mathbf{K}_{c_{1,\varepsilon},\dots,c_{n,\varepsilon}}\varphi,\qquad
    \forall\,\varphi\in\widetilde{\mathbb{H}}{}^r_p(\mathbb{R}^+),
 $$
where
\begin{align} \label{eqn50}
\mathbf{K}_{c_{1,\varepsilon},\dots,c_{n,\varepsilon}}\varphi(t):=&\int_0^\infty
    \mathcal{K}_{c_{1,\varepsilon},\dots,c_{n,\varepsilon}}\left(\frac t\tau\right)\varphi(\tau)\frac{d\tau}\tau=\sum_{j=1}^nd_j(\varepsilon)
    \mathbf{K}^1_{c_{j,\varepsilon}}\varphi(t),\nonumber\\
\mathcal{K}_{c_{1,\varepsilon},\dots,c_{m,\varepsilon}}(t):=&\frac1{(t-c_{1,\varepsilon})
    \cdots(t-c_{n,\varepsilon})}=\sum_{j=1}^n\frac{d_j(\varepsilon)}{t-c_{j,\varepsilon}},\\
c_{j,\varepsilon}\!=\!c(1\!+\!\varepsilon e^{i\omega_j}),& \;
\omega_j\!\in\!(-\pi,\pi), \; \arg
    c_{j,\varepsilon},\ \arg (c_{j,\varepsilon}\,\gamma)\!\not=\!0,\; j\!=\!1,\dots,m. \nonumber
\end{align}
Since $n\in \{3,4,\ldots\}$ the argument $\arg c$ does not vanish.
Hence, the points $\omega_1,\dots,\omega_n\in(-\pi,\pi]$ are
pairwise different, i.e., $\omega_j\not=\omega_k$ for $j\not=k$. By
equating the numerators in the formula \eqref{eqn50} we find the
coefficients $d_1(\varepsilon),\ldots,d_{n-1}(\varepsilon)$.

Note that the operators $\mathbf{K}^3_c,\mathbf{K}^4_c,\ldots$
appear rather rarely in applications. Therefore, in this work exact
formulae are given in the case of the operators $\mathbf{K}^1_c$ and
$\mathbf{K}^2_c$ only.
\end{remark}

Assume that $a_0,\ldots,a_n,b_1,\ldots,b_n\in
C\mathfrak{M}_p(\overline{\mathbb{R}} \setminus\{0\})$,
$c_1,\ldots,c_n\in \mathbb{C}$ and consider the model operator
$\mathbf{A}:
\mathbb{H}^s_p(\mathbb{R}^+)\to\mathbb{H}^s_p(\mathbb{R}^+)$,
\begin{equation}\label{eqn51}
\mathbf{A}:=d_0I+W_{a_0}+\sum_{j=1}^nW_{a_j}\mathbf{K}^1_{c_j}W_{b_j},
\end{equation}
comprising the identity $I$, Fourier $W_{a_0},\ldots,W_{a_n}$,
$W_{b_1},\ldots,W_{b_n}$ and Mellin
$\mathbf{K}^1_{c_1},\ldots,\mathbf{K}^1_{c_n}$ convolution
operators. In order to ensure proper mapping properties of the
operator
$\mathbf{A}\;:\;\widetilde{\mathbb{H}}^s_p(\mathbb{R}^+)\to\mathbb{H}^s_p(\mathbb{R}^+)$,
we additionally assume that if $s\leq1/p-1$ or $s\geq1/p$, then the
functions $a_1(\xi),\ldots,a_n(\xi)$ and $b_1(\xi),\ldots,b_n(\xi)$
have bounded analytic extensions in the lower ${\rm Im}\,\xi<0$ and
the upper  ${\rm Im}\,\xi>0$ half planes, correspondingly.

If $1/p-1<s<1/p$, then the spaces
$\widetilde{\mathbb{H}}^s_p(\mathbb{R}^+)$ and
$\mathbb{H}^s_p(\mathbb{R}^+)$ coincide (can be identified) and the
analytic extendability assumption are not needed. However, we do not
consider this situation here since it requires a special treatment.
Besides, it does not appear in applications.

Now we can describe the symbol $\mathcal{A}^s_p$ of the model
operator $\mathbf{A}$. According to the formulae \eqref{eqn40} and
\eqref{eqn41} one has
\begin{equation}\label{eqn52}
\mathcal{A}^s_p(\omega):=d_0\mathcal{I}^s_p(\omega)+\mathcal{W}^s_{a_0,p}(\omega)
     +\sum_{j=1}^n\mathcal{W}^0_{a_j,p}(\omega)\mathcal{K}^{1,s}_{c_j,p}(\omega)
     \mathcal{W}^0_{b_j,p}(\omega),
\end{equation}
where the symbols $\mathcal{I}^s_p(\omega)$,
$\mathcal{W}^0_{a,p}(\omega)$, $\mathcal{W}^s_{a,p}(\omega)$ and
$\mathcal{K}^{1,s}_{c_j,p}(\omega)$ have the form
\begin{subequations}
 \begin{align}\label{eqn53a}
\mathcal{I}^s_p(\omega)&:=\begin{cases}
    g^s_{-\gamma,\gamma,p}(\infty,\xi), & \omega=(\xi,\infty)\in\overline{\Gamma}_1,
    \\[1ex]
\left(\displaystyle\frac{\eta-\gamma}{\eta+\gamma}\right)^{\mp s}, &
     \omega=(+\infty,\eta)\in\Gamma^\pm_2, \\[1ex]
     e^{\pi si}, &\omega=(\xi,0)\in\overline{\Gamma}_3,\end{cases}\\[2ex]
\label{eqn53bb} \mathcal{W}^0_{a,p}(\omega)&:=\begin{cases}
        a_p(\infty,\xi), & \omega=(\xi,\infty)\in\overline{\Gamma}_1, \\
        a(\mp\eta), & \omega=(+\infty,\eta)\in\Gamma^\pm_2, \\
        a_p(0,\xi), &
        \omega=(\xi,0)\in\overline{\Gamma}_3,\end{cases}\\[2ex]
\label{eqn53b} \mathcal{W}^s_{a,p}(\omega)&:=\begin{cases}
        a^s_p(\infty,\xi), & \omega=(\xi,\infty)\in\overline{\Gamma}_1, \\
        a(\mp\eta)\left(\displaystyle\frac{\eta-\gamma}{\eta+\gamma}\right)^{\mp s}, & \omega=(+\infty,\eta)\in\Gamma^\pm_2, \\
        e^{\pi si}a_p(0,\xi), &
        \omega=(\xi,0)\in\overline{\Gamma}_3,\end{cases}\\[2ex]
 \label{eqn53c}
\mathcal{K}^{1,s}_{c,p}(\omega)&:=\begin{cases}
    \displaystyle\frac{c^{-s}(-c)^{\frac1p-i\xi-1}}{\sin\pi(\frac1p-i\xi)},&\omega=(\xi,\infty)
    \in\overline{\Gamma}_1,,\\[1ex]
    0, &\omega=(\pm\infty,\eta)\in\Gamma^\pm_2, \\[1ex]
    \displaystyle\frac{c^{-s}(-c)^{\frac1p+s-i\xi-1}}{\sin\pi(\frac1p-i\xi)},&\omega=(\xi,0)
    \in\overline{\Gamma}_3,\qquad \text{for} \quad 0<|\arg(c\,\gamma)|<\pi,\end{cases}\\[1.5ex]
a^s_p(\infty,\xi)&:=\frac{e^{2\pi
si}a(\infty)+a(-\infty)}2+\frac{e^{2\pi si}
    a(\infty)-a(-\infty)}{2i}\cot\pi\Big(\frac1p-i\xi\Big),\nonumber
 \end{align}
 \begin{align*}
&\hskip-7mm
a_p(x,\xi):=\frac{a(x+0)+a(x-0)}2+\frac{a(x+0)-a(x-0)}{2i}
     \cot\pi\Big(\frac1p-i\xi\Big),\quad x=0,\infty,\\[1.5ex]
&\hskip-7mm g^s_{-\gamma,\gamma,p}(\infty,\xi):=\frac{e^{2\pi
si}+1}2
     +\frac{e^{2\pi si}-1}{2i}\cot\pi\Big(\frac1p-i\xi\Big)=e^{\pi si}\frac{\sin\pi\Big(\frac1p+s-i\xi\Big)}
     {\sin\pi\Big(\frac1p-i\xi\Big)},\\
&\hskip88mm\xi\in\mathbb{R},\quad \eta\in\mathbb{R}^+,\nonumber
 \end{align*}
where
 \[
-\pi\leqslant\arg\,c<\pi,\quad\arg\,c\not=0,\quad-\pi<\arg(c\,\gamma_0)<0,\quad
0<\arg\gamma,\arg\gamma_0<\pi,
 \]
and $c^s=|c|^se^{is\arg\,c}$, $(-c)^\delta=|c|^\delta
e^{-i\delta\arg\,c}$ for $c,\delta\in\mathbb{C}$.

In the case where $a(-\infty)=1$ and $a(+\infty)=e^{2\pi\alpha i}$
the symbol $a^s_p(\infty,\xi)$ takes the form
 \begin{equation}\label{eqn53f}
a^s_p(\infty,\xi)=e^{\pi(s+\alpha)i}\frac{\sin\pi\left(\frac1p+s+\alpha-i\xi\right)}{
     \cos\pi\Big(\frac1p-i\xi\Big)}.
 \end{equation}
\end{subequations}

Note, that the Mellin convolution operator $\mathbf{K}^1_{-1}$,
\begin{equation*}
\begin{array}{c}
\displaystyle
\mathbf{K}^1_{-1}\varphi(t):=\frac1\pi\int\limits_0^\infty\displaystyle\frac{\varphi(\tau)
     \,d\tau}{t+\tau}=\mathfrak{M}^0_{\mathcal{M}_\frac1p\mathcal{K}^1_{-1}},\quad
\mathcal{M}_\frac1p\mathcal{K}^1_{-1}(\xi)=\displaystyle\frac{1
}{\sin\pi\left(\frac1p-i\xi \right )},
 \end{array}
 \end{equation*}
which often appears in applications, has a rather simple symbol if
considered in the Bessel potential space
$\mathbb{H}^s_p(\mathbb{R}^+)$. Thus using formula \eqref{eqn53c}
with $c=-1$, one obtains
 \begin{align*}
\hskip-5mm\mathcal{K}^{1,s}_{-1,p}(\omega):=\begin{cases}
     \displaystyle\frac{e^{\pi si}}{\sin\pi(\beta-i\xi)},&\omega=(\xi,\infty) \in\overline{\Gamma}_1 \cup\overline{\Gamma}_3,\\
     0, &\omega=(\pm\infty,\eta)\in\Gamma^\pm_2,. \end{cases}
 \end{align*}
 %
\begin{theorem}\label{t5.4}
Let $1<p<\infty$, $s\in\mathbb{R}$. The operator
 \begin{equation}\label{eqn54}
\mathbf{A}:\widetilde{\mathbb{H}}{}^s_p(\mathbb{R}^+)\longrightarrow
     \mathbb{H}^s_p (\mathbb{R}^+)
 \end{equation}
defined in \eqref{eqn51} is Fredholm if and only if its symbol
$\mathcal{A}^s_p(\omega)$ described by the relations \eqref{eqn52},
\eqref{eqn53a}--\eqref{eqn53f}, is elliptic. If $\mathbf{A}$ \, is
Fredholm, then
\begin{equation*}
    {\rm Ind}\mathbf{A}=-{\rm ind}\det\mathcal{A}^s_p.
\end{equation*}
\end{theorem}
\textbf{Proof.} Let $c_j, d_j\in\mathbb{C}$,
$-\pi\leqslant\arg\,c_j<\pi$ $\arg\,c_j\not=0$, for $j=1,\ldots,n$.
Lifting $\mathbf{A}$ up to an operator on the space
$\mathbb{L}_p(\mathbb{R}^+)$ we get
\begin{eqnarray}\label{eqn51a}
\mathbf{\Lambda}^s_{-\gamma}\mathbf{A}\mathbf{\Lambda}^{-s}_{\gamma}
     =d_0\mathbf{\Lambda}^s_{-\gamma}\mathbf{\Lambda}^{-s}_{\gamma}
     +\mathbf{\Lambda}^s_{-\gamma}W_{a_0}\mathbf{\Lambda}^{-s}_{\gamma}
+\sum_{j=1}^nW_{a_j}\mathbf{\Lambda}^s_{-\gamma}\mathbf{K}^1_{c_j}
     \mathbf{\Lambda}^{-s}_{\gamma}W_{b_j},
\end{eqnarray}
where $c^{-s}=|c|^{-s}e^{-is\arg\,c}$ and $\gamma$ is such that
$0<\arg\,\gamma<\pi$, $-\pi<\arg(c_j\,\gamma)<0$ for all
$j=m+1,\ldots,n$.

In \eqref{eqn51a} we used special properties of convolution
operators, namely,
 $$
 \mathbf{\Lambda}^s_{-\gamma}W_{a_j}=W_{a_j}\mathbf{\Lambda}^s_{-\gamma}, \quad
 W_{b_j}\mathbf{\Lambda}^s_\gamma=\mathbf{\Lambda}^s_\gamma
W_{b_j}, \quad \mathbf{\Lambda}^{\mp s}_{\pm\gamma}=W_{\lambda^{\mp
s}_{\pm\gamma}},
 $$
which follows from the analytic extendability of the functions
$\lambda^s_{-\gamma}, a_1(\xi),\ldots,a_n(\xi)$ and
$\lambda^{-s}_\gamma,b_1(\xi),\ldots,b_n(\xi)$ into the lower ${\rm
Im}\,\xi<0$ and upper ${\rm Im}\,\xi>0$ half planes, respectively.

The model operators $I$, $W_{a}$ and $\mathbf{K}^1_c$ lifted to the
space $\mathbb{L}_p(\mathbb{R}^+)$ have the form
 \begin{align}\label{eqn55}
&\mathbf{\Lambda}^s_{\gamma}I\mathbf{\Lambda}^{-s}_{\gamma}=W_{g^s_{-\gamma,\gamma}},
     \qquad \mathbf{\Lambda}^s_{\gamma}W_{a}\mathbf{\Lambda}^{-s}_{\gamma}
     =W_{ag^s_{-\gamma,\gamma}},\nonumber\\[1ex]
&\mathbf{\Lambda}^s_{\gamma}\mathbf{K}^1_c\mathbf{\Lambda}^{-s}_{\gamma}
     =\left\{\begin{array}{ll}c^{-s}\mathbf{K}^1_c W_{g^s_{-c\,\gamma,\gamma}}\hfill{\text{if}}&
     -\pi<\arg(c\,\gamma)<0,\\[2ex]
&\quad 0<\arg(c\,\gamma)<\pi,\\[-1.5ex]
     c^{-s}\mathbf{K}^1_c W_{g^s_{-\,\gamma,-\gamma_0}g^s_{-c\,\gamma_0,\gamma}}+\mathbf{T},\quad \text{if}
     &\\[-1.5ex]
     &\hskip3mm -\pi<\arg(c\,\gamma_0)|<0,
     \end{array}\right.
 \end{align}
where $\mathbf{T}$ is a compact operator. Here, as above,
$-\pi\leqslant\arg\,c<\pi$, $\arg\,c\not=0$, $0<\arg\,\gamma<\pi$,
$0<\arg\,\gamma_0<\pi$ and either $-\pi<\arg(c\,\gamma)<0$ or, if
$-\pi<\arg(c\,\gamma)<0$, then $-\pi<\arg(c\,\gamma_0)|<0$. Recall
that $c^{-s}=|c|^{-s}e^{-is\arg\,c}$.

Therefore, the operator
$\mathbf{\Lambda}^s_{-\gamma}\mathbf{A}\mathbf{\Lambda}^{-s}_{
\gamma}$ in \eqref{eqn51a} can be rewritten as follows
\begin{eqnarray}\label{eqn51b}
\mathbf{\Lambda}^s_{-\gamma}\mathbf{A}\mathbf{\Lambda}^{-s}_{\gamma}
     =d_0W_{g^s_{-\gamma,\gamma}}+W_{a_0g^s_{\gamma,\gamma}}+\sum_{j=1}^mc^{-s}_j
     W_{a_j}\mathbf{K}^1_{c_j}W_{g^s_{-c_j\gamma,-\gamma}}W_{b_j}\nonumber\\
+\sum_{j=m+1}^nc^{-s}_jW_{a_j}\mathbf{K}^1_{c_j}W_{g^s_{-\gamma,-\gamma_j}
     g^s_{-c_j\gamma_j,\gamma}}W_{b_j}+\mathbf{T}\;:\;\mathbb{L}_p(\mathbb{R}^+)
     \longrightarrow\mathbb{L}_p(\mathbb{R}^+),
\end{eqnarray}
where $\mathbf{T}$ is a compact operator and we ignore it when
writing the symbol of $\mathbf{A}$.

Now we define the symbol of the initial operator $\mathbf{A}\;:\;
\widetilde{\mathbb{H}}^s_p(\mathbb{R}^+)\rightarrow
\mathbb{H}^s_p(\mathbb{R}^+)$ of \eqref{eqn51} as the symbol of the
corresponding lifted operator
$\mathbf{\Lambda}^s_{-\gamma}\mathbf{A}\mathbf{\Lambda}^{-s}_{\gamma}:
\mathbb{L}_p(\mathbb{R}^+)\rightarrow\mathbb{L}_p(\mathbb{R}^+)$ of
\eqref{eqn51b}.

To write the symbol of the lifted operator in the Lebesgue space
$\mathbb{L}_p(\mathbb{R}^+)$ let us first find the limits of
involved functions (symbols). The function $g^s_{-\gamma,\gamma}\in
C(\mathbb{R})$ is continuous on $\mathbb{R}$, but has different
limits at the infinity, viz.,
 \begin{align}\label{eqn56a}
 g^s_{-\gamma,\gamma}(-\infty)=1, \quad g^s_{-\gamma,\gamma}(+\infty)=e^{2\pi si},
      \qquad g^s_{-\gamma,\gamma}(0)=e^{\pi si},
 \end{align}
while the functions $g^s_{-\gamma,-\gamma_0},\
g^s_{-c\gamma,\gamma}, g^s_{-c\gamma_0,\gamma} \in C(\mathbb{R})$
are continuous on $\mathbb{R}$ including the infinity. Thus
 \begin{align}\label{eqn56b}
 \begin{array}{c}
g^s_{-c\,\gamma,\gamma}(\pm\infty)=g^s_{-\gamma,-\gamma_0}(\pm\infty)
    =g^s_{-c\,\gamma_0,\gamma}(\pm\infty)=1,\\[2ex]
g^s_{-\gamma,-\gamma_0}(0)g^s_{-c\,\gamma_0,\gamma}(0)
     =\left(\displaystyle\frac{-\gamma}{-\gamma_0}\right)^s
     \left(\displaystyle\frac{-c\gamma_0}{\gamma}\right)^s
     =\left(-c\right)^s,\\[2.5ex]
g^s_{-c\,\gamma,\gamma}(0)=\left(-c\right)^s\quad\text{if}\quad
     -\pi\leqslant\arg\,c<\pi, \quad \arg\,c\not=0.
 \end{array}
 \end{align}

In the Lebesgue space $\mathbb{L}_p(\mathbb{R}^+)$, the symbols of
the first two operators in \eqref{eqn51b}, are written according the
formulae \eqref{eqn40}--\eqref{eqn41} by taking into account the
equalities \eqref{eqn56a} and \eqref{eqn56b}. The symbols of these
operators have, respectively, the form \eqref{eqn53a} and
\eqref{eqn53b}.

For the operators $W_{a_1},\ldots,W_{a_n}$ and
$W_{b_1},\ldots,W_{b_n}$ we can use the formulae
\eqref{eqn40}--\eqref{eqn41} and write their symbols in the form
\eqref{eqn53bb}.

The lifted Mellin convolution operators
 \[
\mathbf{\Lambda}^s_{\gamma}\mathbf{K}^1_{c_j}\mathbf{\Lambda}^{-s}_{\gamma}
     \;:\;\mathbb{L}_p(\mathbb{R}^+)\longrightarrow,\mathbb{L}_p(\mathbb{R}^+)
 \]
comprise both the Fourier convolution operators
$W_{g^s_{-c_j\,\gamma_0,\gamma}}$ and
$W_{g^s_{-\,\gamma,-\gamma_0}g^s_{-c_j\,\gamma_0, \gamma}}$ and the
Mellin convolution operators $\mathbf{K}^1_{c_j}=\mathfrak{M}^0_{
\mathcal{K}^1_{c_j,p}(\xi)}$, with the symbol
$\mathcal{K}^1_{c_j,p}(\xi)
:=\mathcal{M}_{1/p}\mathcal{K}^1_{c_j}(\xi)$ defined in \eqref{eqn9}
and \eqref{eqn10}. The symbol of the operators
$\mathbf{\Lambda}^s_{\gamma}\mathbf{K}^1_{c_j}
\mathbf{\Lambda}^{-s}_{\gamma}$ from \eqref{eqn55} in the Lebesgue
space $\mathbb{L}_p(\mathbb{R}^+)$ is found according formulae
\eqref{eqn40}--\eqref{eqn41}, has the form \eqref{eqn53c} and is
declared the symbol of
$\mathbf{K}^1_{c_j}\;:\;\widetilde{\mathbb{H}}^s_p(\mathbb{R}^+)
\rightarrow\mathbb{H}^s_p(\mathbb{R}^+)$. The symbols of Fourier
convolution factors $W_{g^s_{-c_j\,\gamma_0,\gamma}}$ and
$W_{g^s_{-\,\gamma,-\gamma_0}g^s_{-c_j\,\gamma_0, \gamma}}$, which
contribute the symbol of $\mathbf{K}^1_{c_j}=\mathfrak{M}^0_{
\mathcal{K}^1_{c_j,p}}$ are written again according formulae
\eqref{eqn40}--\eqref{eqn41} by taking into account the equalities
\eqref{eqn56a} and \eqref{eqn56b}.

To the lifted operator applies Theorem \ref{t4.2} and gives the
result formulated in Theorem \ref{t5.4}. \rbx

In the proof of the foregoing Theorem \ref{t5.4}, a local principle
is used. As a byproduct, a result which itself is important in
applications is obtained. We formulate it separately as a corollary.
Note that the definition of the local invertibility and a short
introduction to a local principle can be found in \cite{GK79,Si65}.
 %
\begin{corollary}\label{c5.4}
Let $1<p<\infty$, $s\in\mathbb{R}$. The operator
 \begin{equation*}
\mathbf{A}:\widetilde{\mathbb{H}}{}^s_p(\mathbb{R}^+)\longrightarrow
     \mathbb{H}^s_p (\mathbb{R}^+),
 \end{equation*}
defined in \eqref{eqn51}, is locally invertible at
$0\in\mathbb{R}^+$ if and only if its symbol
$\mathcal{A}^s_p(\omega)$, defined in \eqref{eqn52},
\eqref{eqn53a}--\eqref{eqn53f}, is elliptic on $\Gamma_1$, i.e.
 \[
\inf_{\omega\in\Gamma_1}\left|\det\,\mathcal{A}^s_p(\omega)\right|
     =\inf_{\xi\in\mathbb{R}}\left|\det\,\mathcal{A}^s_p(\xi,\infty)\right|>0.
 \]
\end{corollary}

The next results are concerned with the operators acting in the
Sobolev--Slobodeckij (Besov) spaces. For the definition of the
corresponding spaces
$\mathbb{W}^s_p(\Omega)=\mathbb{B}^s_{p,p}(\Omega)$,
$\widetilde{\mathbb{W}}^s_p(\Omega)=\widetilde{\mathbb{B}}{}^s_{p,p}(\Omega)$
for an arbitrary domain $\Omega\subset\mathbb{R}^n$, including the
semi-axis $\mathbb{R}^+$, we refer the reader to the monograph
\cite{Tr95}.
 %
\begin{corollary}\label{c5.5}
Let $1<p<\infty$, $s\in\mathbb{R}$. If the operator
$\mathbf{A}:\widetilde{\mathbb{H}}{}^s_p(\mathbb{R}^+)
\longrightarrow\mathbb{H}^s_p (\mathbb{R}^+)$, defined in
\eqref{eqn51}, is Fredholm (invertible)  for all $s\in(s_0,s_1)$ and
$p\in(p_0,p_1)$, where $-\infty<s_0<s_1 <\infty$,
$1<p_o<p_1<\infty$, then the operator
\begin{equation}\label{eqn57}
    \mathbf{A}:\widetilde{\mathbb{W}}{}^s_p (\mathbb{R}^+) \longrightarrow \mathbb{W}^s_p (\mathbb{R}^+),\;\;\; s\in(s_0,s_1), \;\; p\in(p_0,p_1)
\end{equation}
is Fredholm (invertible) in the Sobolev--Slobodeckij (Besov) spaces
$\mathbb{W}^s_p =\mathbb{B}^s_{p,p}$, and
\begin{equation}\label{eqn58}
    {\rm Ind}\mathbf{A}=-{\rm ind}\det\mathcal{A}^s_p.
\end{equation}
\end{corollary}

\textbf{Proof.} Recall that the Sobolev--Slobodeckij (Besov) spaces
$\mathbb{W}^s_p=\mathbb{B}^s_{p,p}$ emerge as the result of
interpolation with the real interpolation method between Bessel
potential spaces

\begin{equation}\label{eqn59}
\begin{aligned}
\big(\mathbb{H}_{p_0}^{s_0}(\Omega),\mathbb{H}_{p_1}^{s_1}(\Omega)\big)_{\theta,p}
& =\mathbb{W}_p^s(\Omega),\;\;
    s:=s_0(1-\theta)+s_1\theta, \\
\big(\widetilde{\mathbb{H}}{}_{p_0}^{s_0}(\Omega),\widetilde{\mathbb{H}}{}_{p_1}^{s_1}(\Omega)\big)_{\theta,p}
    & =\widetilde{\mathbb{W}}{}_p^s(\Omega),\;\;p:=\frac1{p_0}\,(1-\theta)+\frac1{p_1}\,\theta, \;\; 0<\theta<1.
\end{aligned}
\end{equation}

If $\mathbf{A}:\widetilde{\mathbb{H}}{}^s_p
(\mathbb{R}^+)\longrightarrow \mathbb{H}^s_p (\mathbb{R}^+)$ is
Fredholm (invertible) for all $s\in(s_0,s_1)$ and $p\in(p_0,p_1)$,
it has a regularizer $\mathbf{R}$ (the inverse
$\mathbf{A}^{-1}=\mathbf{R}$, respectively), which is bounded in the
setting
 \[
\mathbf{R}:\mathbb{W}^s_p (\mathbb{R}^+)
\longrightarrow\widetilde{\mathbb{W}}{}^s_p (\mathbb{R}^+)
 \]
due to the interpolation \eqref{eqn59} and
 \[
\mathbf{R}\mathbf{A}=I+\mathbf{T}_1, \quad
\mathbf{A}\mathbf{R}=I+\mathbf{T}_2,
 \]
where $\mathbf{T}_1$ and $\mathbf{T}_2$ are compact in
$\widetilde{\mathbb{H}}{}^s_p (\mathbb{R}^+)$ and in
$\mathbb{H}^s_p(\mathbb{R}^+)$, or $\mathbf{T}_1=\mathbf{T}_2=0$ if
$\mathbf{A}$ is invertible.

Due to the Krasnoselskij interpolation theorem (see \cite{Tr95}),
$\mathbf{T}_1$ and $\mathbf{T}_2$ are compact in
$\widetilde{\mathbb{W}}{}^s_p(\mathbb{R}^+)$ and in $\mathbb{W}^s_p
(\mathbb{R}^+)$, respectively for all $s\in(s_0,s_1)$ and
$p\in(p_0,p_1)$ and, therefore, $\mathbf{A}$ in \eqref{eqn57} is
Fredholm (is invertible, respectively).

The index formulae \eqref{eqn58} follows from the embedding
properties of the Sobolev--Slobodeckij and Bessel potential spaces
by standard well-known arguments. \rbx

 \end{document}